
\input psfig
\input amssym.def
\input amssym
\magnification=1100
\baselineskip = 0.25truein
\lineskiplimit = 0.01truein
\lineskip = 0.01truein
\vsize = 8.5truein
\voffset = 0.2truein
\parskip = 0.10truein
\parindent = 0.3truein
\settabs 12 \columns
\hsize = 5.4truein
\hoffset = 0.5truein

\setbox\strutbox=\hbox{%
\vrule height .708\baselineskip
depth .292\baselineskip
width 0pt}
\font\caps=cmcsc10

\def\sqr#1#2{{\vcenter{\vbox{\hrule height.#2pt
\hbox{\vrule width.#2pt height#1pt \kern#1pt
\vrule width.#2pt}
\hrule height.#2pt}}}}
\def\square{\mathchoice\sqr46\sqr46\sqr{3.1}6\sqr{2.3}4}
\def\leaderfill{\leaders\hbox to 1em{\hss.\hss}\hfill}
\font\bigtenrm=cmr10 scaled 1400
\tenrm

\vskip 2in
\centerline{\bf {\bigtenrm WORD HYPERBOLIC DEHN SURGERY}}
\tenrm
\vskip 14pt
\centerline{MARC LACKENBY}
\vskip 18pt

\centerline{\caps 1. Introduction}
\vskip 6pt

In the late 1970's, Thurston dramatically changed the nature of 3-manifold
theory with the introduction of his Geometrisation Conjecture, and by
proving it in the case of Haken 3-manifolds [24].
The conjecture for general closed orientable 3-manifolds
remains perhaps the most important unsolved problem in the subject.
A weaker form of the conjecture [20] deals with the fundamental
group of a closed orientable 3-manifold. It 
proposes that either it contains ${\Bbb Z} \oplus
{\Bbb Z}$ as a subgroup or it is word hyperbolic,
in the sense of Gromov [11]. Word hyperbolic groups are
precisely those groups which satisfy a 
linear isoperimetric inequality. They are of fundamental
importance in geometric group theory and have very many
useful properties. 

Many non-Haken 3-manifolds are known to satisfy the geometrisation
conjecture, due to Thurston's hyperbolic Dehn surgery theorem [23].
This asserts that, if one starts with a compact orientable 3-manifold $M$
with $M - \partial M$ supporting a complete finite volume
hyperbolic structure, and one Dehn fills each
component of $\partial M$, then one obtains a hyperbolic
3-manifold, providing that a finite number of `exceptional'
slopes are avoided on each component of $\partial M$.
It remains an interesting unsolved problem to establish
how large this collection of exceptional slopes can be.
Hodgson and Kerckhoff have shown [13] that there is an upper bound
on the number of exceptional slopes, the upper bound being
a (large) number independent of $M$. In this paper, we
will show that very many Dehn fillings yield a 3-manifold
which is irreducible, atoroidal and not Seifert fibred, and has infinite,
word hyperbolic fundamental group. We will approach this problem in a 
number of different ways. 

Our first technique is differential geometric. We will
establish an extension of Thurston and Gromov's $2\pi$ theorem [4],
which we now describe. Pick a horoball neighbourhood $N$
of the cusps of $M - \partial M$. Then, with respect to $N$,
each slope $s$ on $\partial M$ inherits a {\sl length} which is
defined to be the length of the shortest curve on $\partial N$ with
slope $s$. Thurston and Gromov showed that if the Dehn filling
slope on each component of $\partial M$ has length more than
$2 \pi$, then the 3-manifold obtained by Dehn filling has
a negatively curved Riemannian metric. This implies that it is
irreducible, atoroidal and not Seifert fibred, and has infinite, word hyperbolic
fundamental group. The following result improves the critical
slope length from $2\pi$ to $6$.

\noindent {\bf Theorem 3.1.} {\sl Let $M$ be a compact orientable 3-manifold with
interior having a complete finite
volume hyperbolic structure. Let $s_1, \dots, s_n$ be slopes
on $\partial M$, with one $s_i$ on each component of $\partial M$. Suppose that there
is a horoball neighbourhood $N$ of the cusps of $M - \partial M$ on which
each $s_i$ has length more than $6$. Then, the manifold
obtained by Dehn filling along $s_1, \dots, s_n$ is
irreducible, atoroidal and not Seifert fibred, and has infinite, word hyperbolic
fundamental group.}

This result has been established independently by Agol [2] (with
the exception of word hyperbolicity, at this stage). An argument of
Agol [2] gives the following corollary.

\noindent {\bf Corollary 3.4.} {\sl Let $M$ be a compact
orientable 3-manifold with a single torus boundary component,
and with interior supporting a complete finite volume hyperbolic
structure. Then all but at most 12 Dehn fillings on $\partial M$
yield a 3-manifold which is irreducible, atoroidal and not Seifert 
fibred, and has infinite, word hyperbolic fundamental group.}

Our second approach to word hyperbolic Dehn surgery is
more combinatorial. We will define a structure
which we call an `angled ideal triangulation' on $M$,
which was first studied by Casson.
An {\sl angled ideal triangulation} of $M$
is an expression of $M - \partial M$ as a union of ideal 3-simplices
glued along their faces, with each edge of each ideal
3-simplex having an associated interior dihedral angle
in the range $(0, \pi)$.
We insist that the three dihedral angles at each ideal vertex
of each ideal 3-simplex sum to $\pi$, and that the sum of
the dihedral angles around each edge is $2 \pi$.

When $M$ has an angled ideal triangulation, each
boundary component of $M$ inherits
a triangulation, with each triangle having angles
summing to $\pi$. Define the length of
an edge in $\partial M$ to be $\min \lbrace \alpha_1, \dots,
\alpha_6 \rbrace /2$, where
$\alpha_1, \dots, \alpha_6$ are the angles for the
two triangles adjacent to that edge. Each simplicial path in
$\partial M$ then inherits a length. We will
define the `combinatorial length' of a slope on $\partial M$, which will be
{\sl at least} the length of the shortest possible curve
in the 1-skeleton of $\partial M$ having that slope.
Then we will prove (a more general version of) the following theorem, which is
a combinatorial analogue of the $2\pi$ theorem.

\noindent {\bf Theorem 4.9.} {\sl Let $M$ be a
3-manifold with an angled ideal triangulation.
Let $s_1, \dots, s_n$ be a collection of slopes
on $\partial M$, with one $s_i$ on each component of $\partial M$, and each
with combinatorial length more than $2\pi$.
Then the manifold obtained by Dehn filling $M$
along these slopes is irreducible, atoroidal and not Seifert fibred,
and has infinite, word hyperbolic fundamental group.}

\vskip 12pt
\centerline{\psfig{figure=edgelnth.ps}}
\vskip 12pt
\centerline{Figure 1.}

Any complete finite volume non-compact hyperbolic 3-manifold 
can be expressed as a union of convex hyperbolic ideal polyhedra glued
isometrically along their faces [7]. This can usually be decomposed 
further into a hyperbolic ideal triangulation, which determines an 
angled ideal triangulation. However, it is {\sl not} in general true
that an angled ideal triangulation yields a complete
hyperbolic structure. This is
useful, since it implies that Mostow rigidity does not apply to angled
ideal triangulations. One can therefore pick the most convenient
angled ideal triangulation for one's purposes.

Our techniques extend beyond angled ideal triangulations to
expressions of a 3-manifold as a union of ideal polyhedra
with interior angles assigned to each edge. Dually, we will
consider spines with certain extra combinatorial data.
For example, the exterior of an alternating link admits
a well-known spine arising from its diagram.
Using this spine, we will deduce the following result about 
surgery along alternating links.

Let $D$ be an alternating diagram of a knot or link $L$
in $S^3$. We view $D$ as a 4-valent graph $G(D)$ in $S^2$, equipped with
`under-over' crossing information. 
The diagram is {\sl prime} if each simple closed curve in $S^2$ intersecting
$G(D)$ transversely in two points in the interior of edges of $G(D)$ 
divides the 2-sphere into two discs, one of which
contains no crossings of $D$. A region of $D$ is a {\sl bigon region}
if it has precisely two edges of $G(D)$ in its boundary, and an edge of
$G(D)$ is a {\sl bigon edge} if it lies in the boundary of a bigon region.
Define the {\sl twist number} $t(D)$ of $D$ to be half the number
of non-bigon edges of $G(D)$. For example, the alternating
diagram $D$ of the figure-eight knot in Figure 2 has twist number $t(D) = 2$.
More generally, if $K$ is a component of $L$, define
$t(K,D)$ to be half the number of non-bigon edges lying in $K$.
This number may be half-integral, but $t(D)$ is always an integer.
Using this terminology, we can now state one of our main theorems.

\noindent {\bf Theorem 5.1.} {\sl Let $D$ be a connected
prime alternating diagram of a link $L$ in $S^3$. For each
component $K$ of $L$, pick a surgery coefficient $p/q$ (in its lowest terms)
with $\vert q \vert > 8/t(K,D)$. The manifold
obtained by Dehn surgery along $L$, via these surgery
coefficients, is irreducible, atoroidal and not Seifert fibred, and has
infinite, word hyperbolic fundamental group.}

\vskip 12pt
\centerline{\psfig{figure=figeight.ps}}
\vskip 12pt
\centerline{Figure 2.}

Therefore, for `sufficiently complicated' alternating links,
every non-trivial surgery yields a 3-manifold satisfying the conclusions
of the theorem. For example, this is true for all alternating
knots with twist number at least 9. This has the following corollary.

\vfill\eject
\noindent {\bf Corollary 5.2.} {\sl There is a real number $c$ with
the property that if $K$ is an alternating knot whose complement has a
complete hyperbolic structure
with volume at least $c$, then every non-trivial surgery on $K$ yields a
3-manifold which is irreducible, atoroidal and not Seifert fibred, and has
infinite, word hyperbolic fundamental group.}

Surgery along alternating knots has been studied in a number
of different contexts. Menasco and Thistlethwaite [19] embarked
upon an analysis of surfaces with boundary in alternating knot
exteriors, using combinatorial arguments in the spirit of
[18]. They established that, under many circumstances, manifolds
obtained by surgery along alternating knots do not contain
embedded essential 2-spheres or tori. Their bounds on $\vert q \vert$
are a little stronger than those in Theorem 5.1, but their
arguments only seem to work for embedded surfaces, and therefore
do not imply that the fundamental group of the filled-in
3-manifold is infinite and word hyperbolic.
Aitchison, Lumsden and Rubinstein [3] studied certain classes
of alternating links and demonstrated the existence of
non-positively curved metrics on their exteriors. Our
approach is perhaps closest in spirit to theirs, but
there is very little overlap in our results.
Delman [6] and Roberts [22] have shown that every non-trivial surgery along a non-torus 
alternating knot yields a manifold with an essential lamination,
and which is therefore irreducible and has infinite fundamental group [10].
Gabai and Kazez [9] have shown that closed atoroidal 3-manifolds with {\sl genuine} essential
laminations have word hyperbolic fundamental groups, but many
of the laminations constructed by Roberts and Delman are not genuine.

Both our combinatorial approach and our differential geometric
approach to word hyperbolic Dehn surgery have as their basis
Gabai's Ubiquity theorem, which is a
useful tool in showing that a 3-manifold has word
hyperbolic fundamental group. In Section 2,
we will offer a new simplified proof of this theorem.
In Section 3, we will use this to prove Theorem 3.1,
by examining `geodesic spines' of non-compact complete hyperbolic 3-manifolds.
In Section 4, we will examine angled spines
of 3-manifolds (which are generalisations of angled
ideal triangulations) and their relation to word
hyperbolic Dehn surgery. In Section 5, we will apply
these techniques to alternating links. 

\vfill\eject
\centerline{\caps 2. Gabai's Ubiquity theorem}
\vskip 6pt

Recall [11] that a finitely presented group $G$ is
{\sl word hyperbolic} if, for some presentation $\langle
g_1, \dots, g_p \, \vert \, r_1, \dots, r_q \rangle$ of $G$, a linear
isoperimetric inequality is satisfied. In other words, there
exists a constant $c \in {\Bbb R}$,
such that, for every word $w$ in $\langle g_1, \dots, g_p \rangle$
representing the identity element of $G$, we can find an integer 
$n \leq c \vert w \vert$, where $\vert w \vert$ is the length of the word $w$,
and a collection of elements $u_1, \dots, u_n$ of
$\langle g_1, \dots, g_p \rangle$, such that
$$w = \prod_{i=1}^n u_i^{-1} r_{j_i} u_i.$$
Here, the above equality takes place in the free group
$\langle g_1, \dots, g_p \rangle$. It is not hard to see
that if $G$ satisfies a linear isoperimetric
inequality for some presentation, then it will
satisfy such an inequality for all presentations, although
the choice of constant $c$ might differ.
If $M$ is some compact manifold (possibly with boundary)
with a Riemannian metric $g$,
then its fundamental group is word hyperbolic
if and only if, for some constant $c$, each homotopically trivial loop
$K$ in $M$ bounds a disc $D$ with
$${\rm Area}(D) \leq c \ {\rm Length}(K).$$
Again, this property is independent of the
Riemannian metric $g$, but the constant $c$ is not.
Similarly, we can consider a simplicial metric
on $M$, and only consider curves $K$ and discs $D$
which are mapped into $M$ simplicially. Alternatively,
we can consider any metric on $M$ which is
bi-Lipschitz equivalent to a Riemannian metric.

Gabai introduced his Ubiquity theorem in [8].
He and Kazez used it in [9] to show that any closed atoroidal 3-manifold
with a genuine essential lamination has word hyperbolic
fundamental group. Here, we offer a simplified version
of the theorem which is sufficient for our purposes,
together with a new proof using cone manifolds.
This proof has also been observed by Thurston.
We say that a {\sl slope} on a torus is the free homotopy
class of an essential simple closed curve. If $s_1, \dots, s_n$ are
slopes on distinct toral boundary components
of a 3-manifold $M$, then we denote the 3-manifold obtained
by Dehn filling along these slopes by $M(s_1, \dots, s_n)$.

\vfill\eject
\noindent {\bf Theorem 2.1.} {\sl Let $M$ be a
compact orientable 3-manifold whose interior supports a complete finite
volume hyperbolic structure. Pick a Riemannian
(respectively, simplicial) metric $g$ on $M$. Let $s_1, \dots, s_n$ be
slopes on $\partial M$ with one $s_i$ on each component of $\partial M$, and suppose that
the core of each surgery solid torus has infinite order in
$\pi_1(M(s_1, \dots, s_n))$. Suppose that for each smooth (respectively,
simplicial) closed curve $K$ in $M$ which is homotopically trivial in $M(s_1, \dots, s_n)$,
there exists a map of a compact planar surface $F$ into $M$ with
one component of $\partial F$ being sent to $K$ and each
remaining component of $\partial F$ being sent a non-zero multiple of
one of the slopes $s_1, \dots, s_n$, such that
$$\vert F \cap \partial M \vert \leq c \ {\rm Length}(K,g).$$
Here, $c \in {\Bbb R}$ is a constant which may depend
on $M$, $g$ and $s_1, \dots, s_n$, but which is independent of $K$ and $F$.
Then, $\pi_1(M(s_1, \dots, s_n))$ is word hyperbolic.}

The point behind the Ubiquity theorem is that, in order to establish the
word hyperbolicity of $\pi_1(M(s_1, \dots, s_n))$, one need not consider
the area of discs in $M(s_1, \dots, s_n)$, but can instead consider the 
number of times they intersect the cores of the surgery solid tori.

Before we embark on the proof of this theorem, we need to establish some
technical preliminaries. Let $F$ be a compact surface mapped into a 3-manifold
$M$. Then $F$ is {\sl homotopically $\partial$-incompressible} if no
embedded essential arc $R$ in $F$ can be homotoped in $M$ (keeping its
endpoints fixed) to an arc in $\partial M$. Similarly, $F$ is {\sl
homotopically incompressible} if the only simple closed curves
in $F$ which are homotopically trivial in $M$ are those which bound
discs in $F$.

\noindent {\bf Lemma 2.2.} {\sl Let $M$ be a compact orientable
3-manifold with interior supporting a complete finite volume hyperbolic
structure. Let $F$ be a connected compact homotopically $\partial$-incompressible
surface in $M$ with negative Euler characteristic.
Suppose that each component of $\partial F$ either lies in $\partial M$, or
is disjoint from $\partial M$ and not homotopic to a curve in $\partial M$. Suppose
at least one component of $\partial F$ is of the former type, and at most one component
of $\partial F$ is of the latter type. Then 
there is a homotopy of $F - (\partial F \cap \partial M)$ 
to a pleated surface.}

\noindent {\sl Proof.} This is well known, and more details can be
found in [23]. Pick an ideal triangulation of $F - \partial F$, which is
possible since $F$ has negative Euler characteristic and non-empty
boundary. We can ensure that no edge has both endpoints lying in
$\partial F - \partial M$. This determines a lamination on 
$F - (\partial F \cap \partial M)$ by `spinning' the ideal triangulation
of $F - \partial F$ around $\partial F - \partial M$ (if
$\partial F - \partial M \not= \emptyset$) and then adding in
$\partial F - \partial M$. There are (at most) three types
of leaf in this lamination: the curve $\partial F - \partial M$;
edges with both endpoints in $\partial F \cap \partial M$; and
edges with an end spiralling towards $\partial F - \partial M$.
The first of these can be homotoped to a geodesic, since it is not
homotopic to a curve in $\partial M$. The second type of leaf can
homotoped to a geodesic, since $F$ is homotopically 
$\partial$-incompressible. Finally, the third type of leaf can
be homotoped to a geodesic, since the endpoints of a lift of
$\partial F - \partial M$ in ${\Bbb H}^3$ do not lie in the lift of
a cusp of $M - \partial M$. The complement of the lamination is
a collection of ideal triangles. Map each to a totally geodesic ideal
triangle in $M - \partial M$. Then $F  - \partial M$ is a pleated surface.
$\square$

\noindent {\bf Lemma 2.3.} {\sl Let $F$ be a compact orientable incompressible
$\partial$-incompressible surface properly embedded in a compact
orientable 3-manifold $M$ with incompressible boundary. 
Then $F$ is homotopically incompressible and homotopically
$\partial$-incompressible.}

\noindent {\sl Proof.} This is a standard application of the Loop
Theorem [12], which gives us immediately that $F$ is homotopically incompressible.
If $F$ is homotopically $\partial$-compressible, then the doubled surface
$DF$ in the doubled manifold $DM$ is not $\pi_1$-injective, and
so is compressible, by the Loop Theorem. By considering outermost
arcs of intersection between $\partial M$ and a compressing disc
for $DF$, we find that $F$ is either compressible or 
$\partial$-compressible in $M$. $\square$

\noindent {\bf Lemma 2.4.} {\sl Let $F$ be a compact planar
surface mapped into a compact 3-manifold $M$, with $\partial M$
a union of tori. Let $s_1, \dots, s_n$ be slopes on distinct
components of $\partial M$. Suppose that all but one component of $\partial F$ 
are sent to non-zero multiples of $s_1, \dots, s_n$, 
but the remaining component $C$ of $\partial F$
lies on $\partial M$ as a non-zero multiple of a slope other than $s_1, \dots, s_n$.
Then we may find such a surface which is homotopically incompressible
and homotopically $\partial$-incompressible and which has no more
boundary components than $F$.}

\noindent {\sl Proof.} We will show that if $F$ is homotopically
compressible or homotopically $\partial$-compressible, then we
can modify $F$ to a new surface having the properties required
of $F$ but with fewer boundary components. This process will
eventually terminate with the required surface.
Suppose that there is an essential simple
closed curve $R$ in $F$ which is homotopically trivial in $M$. Then
cut $F$ along $R$ and attach a disc to each copy of $R$, yielding
a surface $F'$. The component of $F'$ containing $C$ has
the required properties.

Suppose that there is an embedded essential
arc $R$ in $F$ which can be homotoped (keeping $\partial R$ fixed)
to an arc $R'$ in $\partial M$. Cut $F$ along the arc $R$,
and attach to the new surface two copies of the disc realising the
homotopy between $R$ and $R'$, yielding a surface $F'$.
There are a number of cases to consider.

Suppose first that the two points of $\partial R$ lie in distinct
components $C_1$ and $C_2$ of $\partial F$. Then these components 
are amalgamated into a single boundary component $C'$ of $F'$. 
If neither $C_1$ nor $C_2$ is $C$, then $C'$ is a multiple of one 
of the slopes $s_1, \dots, s_n$. It may be
homotopically trivial, in which case we cap it off with a disc.
If one of $C_1$ and $C_2$ is $C$, then
$C'$ is a non-zero multiple of a slope other 
than $s_1, \dots, s_n$. In each case, we have constructed the
required surface.

Suppose now that the endpoints of $R$ lie in a single component $C_1$ of
$\partial F$. Then the arc $R$ is separating in $F$ (since $F$
is planar). The boundary component $C_1$ becomes two boundary
components $C'_1$ and $C'_2$ of $\partial F'$, which together
are homologous in $\partial M$ to $C_1$. If $C_1 = C$, then
at least one of $C'_1$ and $C'_2$ ($C'_1$, say) is a non-zero multiple of
a slope other than $s_1, \dots, s_n$. Taking the component of $F'$ containing
$C'_1$ gives us the required surface. If $C_1 \not= C$, then
there are two possibilities: either both $C'_1$ and $C'_2$ represent
a multiple of one of the slopes $s_1, \dots, s_n$, or neither do. In the former case, we
take the component of $F'$ containing $C$ (possibly attaching a disc
if one of its boundary components is homotopically trivial). In the latter case,
we take the component of $F'$ not containing $C$. $\square$

Very similar arguments give the following result.

\noindent {\bf Lemma 2.5.} {\sl Let $F$ be a compact planar surface
mapped into a compact 3-manifold $M$, with $\partial M$
a union of tori. Let $s_1, \dots, s_n$ be slopes on distinct
components of $\partial M$, and suppose that the core of
each surgery solid torus in $M(s_1, \dots, s_n)$ has infinite
order in $\pi_1(M(s_1, \dots, s_n))$. Suppose that all but one component of $\partial F$ 
are sent to non-zero multiples of $s_1, \dots, s_n$, but the remaining component of $\partial F$
is sent to a knot $K$ in the interior of $M$. Then we may find such a surface which is 
homotopically incompressible and homotopically $\partial$-incompressible 
and which has no more boundary components than $F$.}

\noindent {\sl Proof of Theorem 2.1.} Let $h$ be the complete finite volume hyperbolic
metric on $M - \partial M$. If the inequality of the theorem holds for
some Riemannian or simplicial metric $g$ on $M$, it holds for every such metric,
after possibly changing our choice of constant $c$.
We are therefore free to give $M$ the Riemannian metric obtained
by removing from $(M - \partial M, h)$ a horoball neighbourhood $N$ of its cusps.
To check the word hyperbolicity of $\pi_1(M(s_1, \dots, s_n))$, it
suffices to verify a linear isoperimetric inequality for
some metric on $M(s_1, \dots, s_n)$ which is bi-Lipschitz equivalent
to a Riemannian metric. We give $M(s_1, \dots, s_n)$ a cone
metric $k$ which is a hyperbolic Riemannian metric away from the
cores $C$ of the surgery solid tori, but which may have
cone angle other than $2 \pi$ at $C$. That such a cone metric
exists on $M(s_1, \dots, s_n)$ is a well known consequence of
the proof of Thurston's hyperbolic Dehn surgery theorem [23].
There is a map from $(M - \partial M,h)$ to $(M(s_1, \dots, s_n), k)$
which does not increase distances by more than a factor of $c_1 > 1$
(say) and is, in fact, bi-Lipschitz with (constant $c_1$) away from $N$.
It collapses a small horoball neighbourhood of the cusps (lying inside $N$)
to the curves $C$.

We now consider a homotopically trivial curve $K$
in $M(s_1, \dots, s_n)$, and will construct a disc bounded by $K$ satisfying
a linear isoperimetric inequality. A homotopy (with area linearly bounded
by ${\rm Length}(K,k)$) pulls $K$ away 
from a regular neighbourhood of $C$. We therefore view $K$ as lying in 
$(M - \partial M) - N$. A further small homotopy takes $K$ 
to a piecewise geodesic curve, each
geodesic segment having length at least $\epsilon$, say,
where $\epsilon$ depends only on $M$.

If $K$ is homotopically trivial in $M - \partial M$, $K$
bounds a disc $D$ in $M - \partial M$ with
${\rm Area}(D,h) \leq c_2 {\rm Length}(K, h)$,
where $c_2 > 0$ depends only on $M$.
This is because $D$ may be realised as a union of
at most $(({\rm Length}(K,h)/ \epsilon) - 2)$ geodesic
triangles, each of which has area at most $\pi$. This implies that
$${\rm Area}(D,k) \leq c_1^2{\rm Area}(D,h)
\leq c_1^2 c_2 {\rm Length}(K,h) \leq c_1^3 c_2 {\rm Length}(K,k),$$
which is the required linear isoperimetric inequality.

Suppose therefore that $K$ is homotopically non-trivial in $M$. As
in the statement of the theorem, there exists a map of a compact
planar surface $F$ into $M$, satisfying the inequality
$\vert F \cap \partial M \vert \leq c \ {\rm Length}(K,h)$.
We may take such an $F$ for which $\vert F \cap \partial M \vert$
is minimal, and so by Lemma 2.5, we may assume that $F$ is homotopically
$\partial$-incompressible. If $\vert F \cap \partial M \vert = 1$, then, as above, we
let $F - (\partial F \cap \partial M)$ be a union of geodesic triangles
and extend it to a disc $D$ in $M(s_1, \dots, s_n)$ with ${\rm Area}(D,k)
\leq c_1^3 c_2 {\rm Length}(K,k)$. If $K$ is homotopic in $M$
to a curve in $\partial M$, then this curve is either a multiple of
some $s_i$ (in which case, we may take $\vert F \cap \partial M \vert = 1$)
or is not a multiple of some $s_i$ (in which case, a core of one of the surgery
solid tori has finite order in $\pi_1(M(s_1, \dots, s_n))$,
contrary to assumption). Hence, we may assume that $F$
satisfies the requirements of Lemma 2.2 and so we may homotope
$F - (\partial F \cap \partial M)$ to a pleated surface $F_1$, say.
This homotopy will take $K$ to a geodesic $K_1$.
It is well known [11] that the annulus $A$ realising the
free homotopy between $K$ and $K_1$ can be taken to have
$${\rm Area}(A,h) \leq c_3 \ {\rm Length}(K, h),$$
for some constant $c_3 > 0$ which depends only on $M$
(by an argument similar to the case where $K$ is homotopically trivial
in $M$). The pleated surface $F_1$ has
$${\rm Area}(F_1, h) = -2 \pi \chi(F_1)
= 2\pi(\vert F \cap \partial M \vert - 1)
< 2\pi c \ {\rm Length}(K,h).$$ 
Therefore, by gluing $F_1$ and $A$, we get that $K$ is 
part of the boundary of a surface $F_2$ as in
the statement of the theorem with 
$${\rm Area}(F_2, h) \leq (2 \pi c + c_3){\rm Length}(K,h).$$ 
This surface $F_2$ extends to a disc $D$ in $M(s_1, \dots, s_n)$, with
$${\rm Area}(D,k) \leq c_1^2(2 \pi c + c_3){\rm Length}(K,h)
\leq c_1^3(2 \pi c + c_3){\rm Length}(K,k),$$
which is the required isoperimetric inequality. $\square$

\vfill\eject
\centerline{\caps 3. An extension of the $2 \pi$ theorem}
\vskip 6pt

The goal of this section is to prove the following result.

\noindent {\bf Theorem 3.1.} {\sl Let $M$ be a compact orientable 3-manifold 
with interior having a complete finite 
volume hyperbolic structure. Let $s_1, \dots, s_n$ be slopes
on $\partial M$, with one $s_i$ on each component of $\partial M$. Suppose that there
is a horoball neighbourhood $N$ of the cusps of $M - \partial M$ on which
each $s_i$ has length more than $6$. Then, the manifold
obtained by Dehn filling along $s_1, \dots, s_n$ is
irreducible, atoroidal and not Seifert fibred, and has infinite, word hyperbolic
fundamental group.}

Our first step is to construct a `geodesic spine' of $M$.
Recall that a {\sl spine} of $M$ is an embedded 2-dimensional cell
complex $S$ with the property that $M$ is a regular neighbourhood of $S$.
A {\sl geodesic spine} is defined to be a spine of $M$ in which
each cell is totally geodesic. Define the {\sl geodesic
spine arising from the horoball neighbourhood $N$}
to be
$$S = \lbrace x \in (M - \partial M) - N : \hbox {$x$ does not have a unique closest
point in $N$} \rbrace.$$
To see that $S$ is a geodesic spine, consider its inverse image $\tilde S$
under the covering map ${\Bbb H}^3 \rightarrow M - \partial M$. In this
cover, $N$ lifts to a collection of disjoint horoballs $\tilde N$.
Let $\tilde E$ be the closure of a component of ${\Bbb H}^3 - \tilde S$.
Each point in $\tilde E - \tilde S$
is closest to a single component of $\tilde N$. Placing this component
of $\tilde N$ at $\infty$ in the upper half space model of 
${\Bbb H}^3$, we see that there is a natural vertical retraction 
$\tilde E \rightarrow \partial \tilde E$.
This retraction is defined for all components of ${\Bbb H}^3
- \tilde S$ and is equivariant with respect to the covering
${\Bbb H}^3 \rightarrow M - \partial M$. Hence, it descends to a retraction
$M - \partial M \rightarrow S$. 
Consider a point of $\tilde S$. If $\lbrace n_1, \dots, n_t \rbrace$
are its nearest points in $\tilde N$, then $n_1, \dots, n_t$ lie
in distinct components $\tilde N_1, \dots, \tilde N_t$ of
$\tilde N$. Also contained in $\tilde S$ are the set of points in ${\Bbb H}^3$
which are equidistant to $\tilde N_1, \dots, \tilde N_t$ and not closer
to any other component of $\tilde N$. This forms a totally geodesic
convex subset of ${\Bbb H}^3$, which we take to be a single cell
in our cell structure on $\tilde S$. This descends to a cell
structure on $S$, which makes it into a geodesic spine.

Note that when $M$ has a single boundary component, 
the geodesic spine is independent of the choice of
horoball neighbourhood $N$ of the cusp. However, when $\partial M$ has
more than one component, the geodesic spine need not be
unique. Possibly the simplest example of this construction is the
geodesic spine $S$ for the figure-eight knot complement, shown in
Figure 3. A cross-section of $\tilde S$ in the upper half space
model of ${\Bbb H}^3$ is shown in Figure 4.

\vskip 24pt
\centerline{\psfig{figure=geodspin.ps,width=4.5in}}
\vskip 12pt
\centerline{Figure 3.}

\vskip 12pt
\centerline{\psfig{figure=extremal.ps}}
\vskip 12pt
\centerline{Figure 4.}

An examination of $S$ in this case will be instructive,
since it satisfies a number of `extremal' properties, 
which will be central to our arguments. Arrange $\tilde N$
so that it touches $\infty$ in the upper half space model for
${\Bbb H}^3$, and arrange $\tilde S$ so that it just touches
$\lbrace (x,y,z) \in {\Bbb H}^3: z = 1 \rbrace$ (which we abbreviate
to $\lbrace z = 1 \rbrace$). Let $z_0$ be a real number in the
range $(\sqrt 3/2,1)$. Note that $\tilde S \cap \lbrace z = z_0 \rbrace$ 
is a collection of disjoint circles. Let $R(z_0)$ be the radius of each circle,
and let $D(z_0)$ be the minimal distance between circles,
measured in the Euclidean metric on $\lbrace z = z_0 \rbrace$.
Note that $2R(z_0) + D(z_0)$ is the distance between
the centres of adjacent circles, which is $1/z_0$.
The numbers $R(z_0)$ and $D(z_0)$ are extremal in the following sense.

\noindent {\bf Lemma 3.2.} {\sl Let $M$ be a compact orientable 3-manifold
whose interior supports a complete finite volume
hyperbolic structure. Let $N$ be a horoball neighbourhood
of its cusps, and let $S$ be the associated geodesic spine.
Let $\tilde N$ and $\tilde S$ be the inverse images of $N$ and $S$ in
${\Bbb H}^3$. Arrange one component $\tilde N_0$ of $\tilde N$ as
$\lbrace z \geq 1 \rbrace$ in the upper half
space model of ${\Bbb H}^3$, and let $\tilde E$ be the closure of the component of
${\Bbb H}^3 - \tilde S$ containing $\tilde N_0$. Let $z_0$ be a real number
in the range $(\sqrt 3 /2, 1)$. Then the intersection $({\Bbb H}^3 - \tilde E) \cap
\lbrace  z = z_0 \rbrace$ is a (possibly empty)
collection of disjoint discs. The radius of each disc
is at most $R(z_0)$, and the distance between
two discs is at least $D(z_0)$.}

\noindent {\sl Proof.} Construct $\tilde E$ as follows.
For each component $\tilde N_i$ of $\tilde N$ other than $\tilde N_0$, let $P_i$
be the totally geodesic plane equidistant between $\tilde N_0$
and $\tilde N_i$. Then $\tilde E$ is the set of all points
above $\bigcup_i P_i$ in the upper half space model.
We will take two such planes $P_1$ and $P_2$ and will
move them in ${\Bbb H}^3$ using the sequence of operations in Figure 5
until they are the planes
$P_1$ and $P_2$ of Figure 4. Each operation will not
decrease the radius of $P_i \cap \lbrace 
z = z_0 \rbrace$, nor will it increase the distance
between $P_1 \cap \lbrace 
z = z_0 \rbrace$ and $P_2 \cap \lbrace 
z = z_0 \rbrace$. Since these circles end up as being
disjoint, they must have started disjoint. Thus, $({\Bbb H}^3 - \tilde E) \cap
\lbrace  z = z_0 \rbrace$ must have started as a (possibly empty)
collection of disjoint discs, satisfying the required properties.

The first operation is to translate $\tilde N_1$ and $\tilde N_2$ 
horizontally until they
touch. The second operation is to scale (in the Euclidean metric on the
upper half space model) both $\tilde N_1$ and
$\tilde N_2$ by the same factor, keeping them just touching one another, 
until at least one of them ($\tilde N_2$, say) also touches $\tilde N_0$. 
We can perform this scale so that the set of points above the old 
$\tilde N_1 \cup \tilde N_2$ contains the
set of points above the new $\tilde N_1 \cup \tilde N_2$. This implies
that the old $P_1 \cup P_2$ lies below the new $P_1 \cup P_2$,
since $P_i$ is the set of points equidistant between $\tilde N_0$
and $\tilde N_i$. Hence,
we have not decreased the radii of $P_i \cap \lbrace 
z = z_0 \rbrace$, nor have we moved $P_1 \cap \lbrace 
z = z_0 \rbrace$ and $P_2 \cap \lbrace 
z = z_0 \rbrace$ further from each other. This
is also true of the third operation, which leaves $\tilde N_2$
unchanged but expands and translates $\tilde N_1$, keeping $\tilde N_2$ and
$\tilde N_1$ just touching each other. We end with $\tilde N_1$
and $\tilde N_2$ as in Figure 5 and hence $P_1$ and $P_2$ as in Figure 4. $\square$

\vskip 12pt
\centerline{\psfig{figure=horomove.ps,width=4in}}
\vskip 12pt
\centerline{Figure 5.}

The basis behind Theorem 3.1 is the following area estimate.
The factor $(\pi/3)$ in this estimate is the reason for the improvement
of the critical slope length from $2 \pi$ to $6$.

\noindent {\bf Lemma 3.3.} {\sl Let $M$ be a compact orientable 3-manifold whose
interior supports a complete finite volume hyperbolic structure.
Let $S$ be a geodesic spine arising from a horoball neighbourhood
$N$ of the cusps of $M - \partial M$. Let $G$ be a compact orientable
(possibly non-embedded) surface in the closure of a component
of $M - S$, with $\partial G \subset \partial M \cup S$ and
with $\partial G \cap \partial M$ a single curve representing
$\pm k [s] \in H_1(\partial M)$, where $k \in {\Bbb N}$ and $s$ is
some slope. Then
$${\rm Area}(G - \partial M) \geq k (\pi/3) \ {\rm Length}(s),$$
where the slope length of $s$ is measured with respect to $N$.}

\noindent {\sl Proof.} Let $E$ be the closure of the component of $(M -
\partial M) - S$ which contains $G - \partial M$, and
let $\tilde E$ be a component of the inverse
images of $E$ in ${\Bbb H}^3$. As in Lemma 3.2, arrange the component $\tilde N_0$
of $\tilde N$ lying in $\tilde E$ as $\lbrace z \geq 1 \rbrace$. 
For $z_0 \geq \sqrt 3 /2$, let $E(z_0)$ be the image of $\tilde E(z_0) = \tilde E
\cap \lbrace  z \geq z_0 \rbrace$ under the covering map ${\Bbb H}^3 
\rightarrow M - \partial M$. We may assume (after a very small
homotopy of $G$) that, for all but finitely many $z_0$, 
$\partial E(z_0) \cap G$ is a collection $C$ of immersed curves in
$\partial E(z_0)$ which is homologous to $\pm k [s]$ 
in $\partial M$. 
For each component of $C$, pick a lift to an arc
or closed curve in $\partial \tilde E(z_0)$, each arc having endpoints on $\tilde S$,
and let $\tilde C'$ be the union of these lifts over all components of $C$. 
Let $\tilde C$ be the image of $\tilde C'$ under the
vertical projection map $\partial \tilde E(z_0) \rightarrow \lbrace  z = z_0
\rbrace$. Discard any components of
$\tilde C$ disjoint from $\tilde E$. Straighten $\tilde C$
inside each disc component of
$({\Bbb H}^3 - \tilde E) \cap \lbrace  z = z_0 \rbrace$, creating
arcs and circles $\tilde C(z_0)$, say,
with total length at least $k \,
{\rm Length}(s)/z_0$. (See Figure 6.)

\vskip 12pt
\centerline{\psfig{figure=circlez0.ps}}
\vskip 12pt
\centerline{Figure 6.}

For our purposes, what is relevant
is the ratio of ${\rm Length}(\tilde C(z_0) \cap \tilde E)$ to 
${\rm Length}(\tilde C(z_0) - \tilde E)$. Suppose initially that
no component of $\tilde C(z_0)$ lies entirely in $\tilde E$.
The arcs of $\tilde C(z_0) \cap \tilde E$
alternate with the arcs of $\tilde C(z_0) - \tilde E$. By Lemma 3.2, each
of the former arcs has length at least $D(z_0)$, and each of
the latter arcs has length at most $2R(z_0)$. Hence,
$${\rm Length}(\tilde C(z_0) \cap \tilde E) \geq {D(z_0) \over D(z_0) + 2R(z_0)}
{\rm Length}(\tilde C(z_0))
\geq D(z_0) k \, {\rm Length}(s).$$
Note that this is true even if components of $\tilde C(z_0)$ lie entirely in 
$\tilde E$. It also remains true for $z_0 \geq 1$ if we define
$D(z_0) = 1/z_0$ and $R(z_0) = 0$. Hence
$$\eqalign{{\rm Area}(G - \partial M) 
&\geq \int_{\sqrt 3 /2}^\infty {{\rm Length}(\tilde C(z) \cap \tilde E) \over z} \, dz \cr
&\geq k \, {\rm Length}(s) \int_{\sqrt 3 /2}^\infty {D(z) \over z} \, dz.\cr}$$
But, it is clear from Figure 7, that the final integral is
precisely one third of the area of an ideal triangle. Therefore
${\rm Area}(G - \partial M) \geq k (\pi/3) \ {\rm Length}(s).$ $\square$

\vfill\eject
\centerline{\psfig{figure=piover3.ps}}
\vskip 12pt
\centerline{Figure 7.}

\noindent {\sl Proof of Theorem 3.1.} Note first that
if $M(s_1, \dots, s_n)$ is irreducible and has infinite, word hyperbolic fundamental group,
then it is neither toroidal nor Seifert fibred.
Suppose that $M(s_1, \dots, s_n)$ is reducible.
Then we can find a compact, incompressible, $\partial$-incompressible,
planar surface $F$ properly embedded in
$M$, with each component of $\partial F$ having one of the
slopes $s_1, \dots, s_n$, and with $\vert \partial F \vert \geq 3$.
By Lemmas 2.3 and 2.2, we may homotope $F - \partial F$ to a pleated 
surface in $M - \partial M$. It then inherits a hyperbolic Riemannian metric
which has 
$${\rm Area}(F - \partial F) = -2 \pi \chi(F) = 2 \pi (\vert F \cap \partial M
\vert - 2) < 2\pi \vert F \cap \partial M \vert.$$
Each component of intersection
between $F$ and $M - S$ intersects $\partial M$ at most once,
otherwise $F$ would be $\partial$-compressible in $M$. Hence, Lemma 3.3
gives that
$${\rm Area}(F - \partial F) \geq \vert F \cap \partial M \vert
(\pi/3) \ \min_{1 \leq i \leq n} {\rm Length}(s_i) > 2 \pi
\vert F \cap \partial M \vert,$$
which gives a contradiction. 

Similarly, we claim that the core of each surgery solid
torus in $M(s_1, \dots, s_n)$ has infinite order in 
$\pi_1(M(s_1, \dots, s_n))$. In particular, 
$\pi_1(M(s_1, \dots, s_n))$ is infinite. For, if not,
some non-zero multiple of some core curve forms the boundary of
a disc mapped into $M(s_1, \dots, s_n)$. We can ensure that
the intersection of this disc with $M$ is a compact planar 
surface $F$ mapped into $M$ with all but one component of $\partial F$ having slope
a non-zero multiple of some $s_i$, and the remaining component
sent to a non-zero multiple of a slope other than $s_1, \dots, s_n$.
By Lemma 2.4, we may take
$F$ to be homotopically $\partial$-incompressible, 
and so we may homotope $F - \partial F$ into pleated form.
The above argument gives us a contradiction. 

In order to prove that $\pi_1(M(s_1, \dots, s_n))$ is word hyperbolic,
we will use Gabai's Ubiquity theorem. We need to give
$M$ some metric; this will be the hyperbolic Riemannian metric on
$(M - \partial M) - {\rm int}(N)$. Consider a curve $K$ in 
$(M - \partial M) - {\rm int}(N)$ which is homotopically trivial
in $M(s_1, \dots, s_n)$. There exists a compact planar surface $F$
in $M$, with $\partial F$ consisting of $K$ and
curves on $\partial M$, each representing a non-zero multiple of 
one of the slopes $s_1, \dots, s_n$.
Since each surgery solid torus in $M(s_1, \dots, s_n)$
has infinite order in $\pi_1(M(s_1, \dots, s_n))$, Lemma 2.5
implies that we may take $F$ to be homotopically 
$\partial$-incompressible. We will show that there is a
constant $c$ and a choice of $F$ with $\vert F \cap \partial M \vert
\leq c \, {\rm Length}(K)$. If $K$ is homotopically trivial in $M$,
then this is immediate. If $K$ is homotopically non-trivial, then there
is a lower bound on its length. Hence, by taking $c$ sufficiently
large, we may assume that $\vert F \cap \partial M \vert \geq 2$.
Hence, by Lemma 2.2, we may
homotope $F - \partial M$ in $M - \partial M$ to a pleated surface $F_1$, taking
$K$ to geodesic $K_1$. Since $F_1$ is pleated, its area is
$${\rm Area}(F_1) = - 2 \pi \chi(F_1)
< 2 \pi \vert F \cap \partial M \vert.$$
As in the proof of Theorem 2.1, we may ensure that the
annulus realising the free homotopy between $K$ and $K_1$ has area
at most $c_3 {\rm Length}(K)$, where $c_3 >0$ is a constant depending
only on $M$. Gluing this annulus to $F_1$ gives a surface
$F_2$ with $${\rm Area}(F_2) < 2 \pi \vert F \cap \partial M \vert +
c_3 {\rm Length}(K).$$

Let $E$ be the closure of some component of $(M - \partial M) - S$.
For any $z_0 \geq \sqrt3/2$, define $E(z_0)$ as in the proof of
Lemma 3.2. Then, there is a vertical projection $E(\sqrt 3 /2) \rightarrow 
\partial E(\sqrt 3/2)$. For each component 
$A$ of $K \cap E(\sqrt 3 /2)$, let $D(A)$ be the vertical
disc in $E(\sqrt 3/2)$ lying below it.
Since $K$ is disjoint from ${\rm int}(N)$, 
$A$ lies in the region $\lbrace \sqrt
3/2 \leq z \leq 1 \rbrace$, and so 
$${\rm Area}(D(A)) \leq (2/ \sqrt 3 -1) {\rm Length}(A).$$

\vfill\eject
\centerline{\psfig{figure=esqrt3.ps}}
\vskip 12pt
\centerline{Figure 8.}

Let $G$ be a component of $F_2 - S$ which touches $\partial M$.
Since $F$ is homotopically $\partial$-incompressible, $G$ must have a single
boundary component which represents $\pm k [s_j] \in H_1(\partial M)$,
for some $s_j \in \lbrace s_1, \dots, s_n \rbrace$ and $k \in {\Bbb N}$.
If $G$ is not disjoint from $K$, let 
$A_1, \dots, A_t$ be the arcs of $K \cap G$. Attach $D(A_1)
\cup \dots \cup D(A_t)$ to $G$ to form a surface $G_1$. (If $G$ is disjoint
from $K$, let $G_1 = G$.) By Lemma 3.3,
$${\rm Area}(G_1) \geq k (\pi/3) \ {\rm Length}(s_j) \geq (\pi/3) \ {\rm Length}(s_j),$$
and so 
$${\rm Area}(G) \geq (\pi/3) \ {\rm Length}(s_j)
- (2/\sqrt 3 -1) {\rm Length}(K \cap G).$$
Summing these inequalities over each component $G$ of $F_2 - S$
which touches $\partial M$, we get that
$${\rm Area}(F_2) \geq \vert F \cap \partial M \vert 
(\pi/3) \min_{1 \leq j \leq n} {\rm Length}(s_j) 
- (2/ \sqrt 3 -1) {\rm Length}(K).$$
Therefore,
$$\eqalign {& 2 \pi \vert F \cap \partial M \vert +
c_3 {\rm Length}(K) > {\rm Area}(F_2)\cr  
&\geq \vert F \cap \partial M \vert 
(\pi/3) \min_{1 \leq j \leq n} {\rm Length}(s_j) 
- (2 / \sqrt 3 -1) {\rm Length}(K),\cr}$$
which gives that
$$\vert F \cap \partial M \vert \Big[
\Big( (\pi/3) \ \min_{1 \leq j \leq n} {\rm Length}(s_j) \Big) - 2 \pi \Big]
\leq (2 / \sqrt 3 -1 + c_3) {\rm Length}(K).$$
Theorem 2.1 now gives that $\pi_1(M(s_1, \dots, s_n))$
is word hyperbolic. $\square$

\noindent {\bf Corollary 3.4.} {\sl Let $M$ be a compact
orientable 3-manifold with a single torus boundary component,
and with interior supporting a complete finite volume hyperbolic
structure. Then all but at most 12 Dehn fillings on $\partial M$
yield a 3-manifold which is irreducible, atoroidal and not Seifert 
fibred, and has infinite, word hyperbolic fundamental group.}

{\sl Proof.} It suffices to show that at most 12 slopes on $\partial M$
have length no more than 6 (with respect to the maximal horoball neighbourhood $N$ of
$\partial M$). This was established by Agol [2], using
a new result of Cao and Meyerhoff [5]. Here is a brief outline
of Agol's argument. He showed first that if two slopes
$s_1$ and $s_2$ on $\partial M$ have length at most six, 
then their distance
$\Delta(s_1, s_2)$  on $\partial M$ is at most 10. This is because
a simple and well-known argument in Euclidean geometry gives that
$$\Delta(s_1, s_2) \leq {{\rm Length}(s_1) \, {\rm Length}(s_2) \over 
{\rm Area}(\partial N)},$$
and Cao and Meyerhoff have recently shown that ${\rm Area}(\partial N)
\geq 3.35$. 
This gives that $\Delta(s_1, s_2) \leq 10$.
Agol then showed that if ${\cal S}$ is a collection of distinct slopes
on a torus, with any two elements of ${\cal S}$ having
distance at most 10, then $| {\cal S} |
\leq 12.$ $\square$

It is worth noting that, in many cases, the bounds of Theorem 3.1
will be far from optimal. The basis behind Theorem 3.1 is the
area estimate in Lemma 3.3, which exploits the fact that the
area of the surface $G$ embedded in the complement of the geodesic spine $S$
picks up more area than the parts of $G$ lying in the horoball
neighbourhood $N$. If the manifold $M - \partial M$ has
large volume in comparison with that of $N$, then there will be parts of $S$
lying far from $N$. If $G$ approaches these parts of $S$, then
its area will be greater than Lemma 3.3 predicts.

\vfill\eject
\centerline{\caps 4. Angled spines of 3-manifolds}
\vskip 6pt

Casson has realised that useful results
may deduced about non-compact finite volume complete hyperbolic 3-manifolds
merely by examining the interior dihedral angles of
a straight ideal triangulation. He studied closed normal
surfaces in these ideal triangulations, and he showed that
normal 2-spheres do not occur, and that the only normal
tori are links of the ideal vertices. This section is
devoted to showing that Casson's ideas may be extended
to surfaces with boundary, and in this way, we
will deduce some interesting results about Dehn
surgery.

Recall that an {\sl ideal triangulation}
of a compact 3-manifold $M$ with non-empty boundary is
an expression of $M - \partial M$ as a union of 3-simplices
glued along their faces, with the vertices then removed.
An {\sl angled ideal triangulation} is an ideal triangulation,
together with the assignment of
a real number in the range $(0, \pi)$ to each
edge of each 3-simplex, known as the {\sl interior angle}
of that edge. We insist that the interior angles around
each edge sum to $2 \pi$, and that the sum of the three
interior angles at each vertex of each 3-simplex is $\pi$.

\vskip 12pt
\centerline{\psfig{figure=tetangle.ps}}
\vskip 12pt
\centerline{Figure 9.}

The dual picture of an ideal triangulation is a special
spine. If the ideal triangulation is angled, then the
associated special spine inherits a certain combinatorial structure.
We generalise this below, by considering spines more general 
than special spines.
A spine of a 3-manifold $M$ with non-empty boundary can
be thickened to a handle decomposition ${\cal H}$ of $M$. 
Hereafter, we will only consider handle structures
arising in this way. We will denote the $i$-handles
of ${\cal H}$ by ${\cal H}^i$. Particularly important will be
the surface ${\cal H}^0 \cap ({\cal H}^1 \cup
{\cal H}^2)$, which we will denote by ${\cal F}$.
Providing each 2-handle of ${\cal H}$ touches some 1-handle,
this surface ${\cal F}$ inherits a handle
structure, with the 0-handles ${\cal F}^0$ of ${\cal F}$ being
${\cal H}^0 \cap {\cal H}^1$, and the 1-handles ${\cal F}^1$ of
${\cal F}$ being ${\cal H}^0 \cap {\cal H}^2$.
We will insist that $H_0 \cap {\cal F}$ is connected
for each 0-handle $H_0$ of ${\cal H}$.
When ${\cal H}$ is dual to an angled ideal triangulation, then each
component of ${\cal F}$ is of the form shown in Figure 10 with
`interior' angles assigned to each 1-handle of ${\cal F}$.
It is this structure which we wish to generalise.

\vskip 12pt
\centerline{\psfig{figure=fangled.ps}}
\vskip 12pt
\centerline{Figure 10.}

There is a well-established theory [15] of normal embedded
surfaces in a handle structure ${\cal H}$ of a 3-manifold $M$, which we now
outline. Any incompressible surface $F$ properly embedded in an irreducible
3-manifold $M$ can be ambient isotoped
so that afterwards, it intersects each handle of ${\cal H}$
in a collection of disjoint discs which respect the product structure
on that handle, as in Figure 11. This gives the surface a handle
structure, with the $i$-handles of $M$ giving rise to the
$i$-handles of $F$.

\vskip 18pt
\centerline{\psfig{figure=stanhand.ps,width=4.5in}}
\vskip 12pt
\centerline{Figure 11.}
\vfill\eject

Furthermore, if $F$ is
$\partial$-incompressible and no component of $F$ is a
disc parallel to a disc in $\partial M$, then we can ensure that
each curve $N$ of $F \cap \partial {\cal H}^0$ is {\sl normal}, 
which means that it is a simple closed curve satisfying each 
of the following conditions:
\item{(i)} $N$ respects the product structures of the 1-handles
of ${\cal F}$;
\item{(ii)} $N$ does not lie entirely in $\partial M$ or
entirely in ${\cal F}^0$;
\item{(iii)} no arc of intersection between $N$ and ${\cal F}^0$
has endpoints lying in the same component of
$\partial {\cal F}^0 \cap {\cal F}^1$, or in the same component
of $\partial {\cal F}^0 - {\cal F}^1$, or in adjacent
components of $\partial {\cal F}^0 \cap {\cal F}^1$
and $\partial {\cal F}^0 - {\cal F}^1$;
\item{(iv)} no arc of intersection between $N$ and $\partial {\cal H}^0
- {\cal F}$ has endpoints lying in the same component of 
$\partial {\cal F}^0 - {\cal F}^1$;
\item{(v)} $N$ intersects any given 1-handle of ${\cal F}$
in at most one arc;
\item{(vi)} $N$ intersects any given component of 
$\partial {\cal H}^0 - {\cal F}$ in at most one arc.

\noindent We will say that a surface $F$ satisfying all of
the above restrictions is {\sl normal}.

\vskip 18pt
\centerline{\psfig{figure=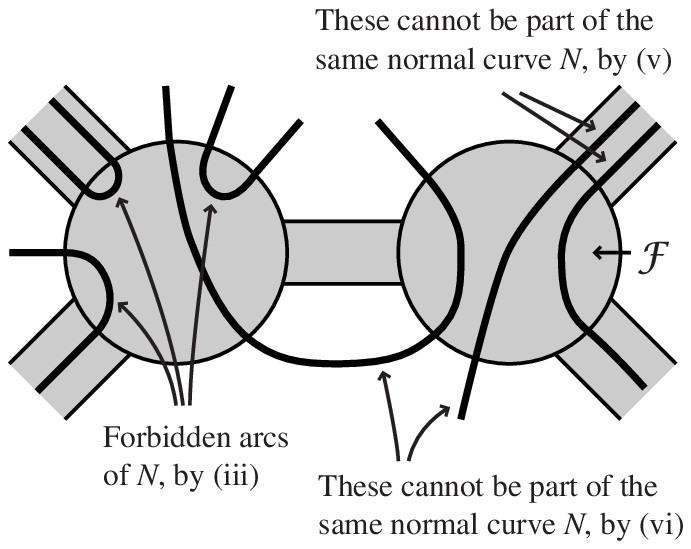}}
\nobreak
\vskip 6pt
\centerline{Figure 12.}

We can now frame some crucial definitions. Consider a spine of
a 3-manifold $M$ with non-empty boundary. Assign
two real numbers to each 1-handle of ${\cal F}$, each
number lying in the range $(0, \pi)$, and the two numbers summing to
$\pi$. These numbers are known as the {\sl interior} and
{\sl exterior angles} of that 1-handle. If $N$ is
a normal curve, suppose that $\epsilon_1, \dots, \epsilon_n$
are the exterior angles of the 1-handles of ${\cal F}$
along which $N$ runs. Then define the {\sl combinatorial area} of $N$ to be
$$a(N) = (\sum_{i=1}^n \epsilon_i) -2 \pi +
\pi \vert N \cap \partial M \vert.$$
\noindent We will say that the angles on ${\cal F}^1$
determine an {\sl angled spine} providing that
\item{$\bullet$} every normal curve in $\partial {\cal H}^0$ 
has non-negative combinatorial area, and
\item{$\bullet$} the interior angles around
any 2-handle of ${\cal H}$ sum to $2 \pi$.

\noindent We define the {\sl combinatorial area} $a(H)$ of a 0-handle $H$
of a normal surface $F$ to be the combinatorial area of $\partial H$.
We also define the {\sl combinatorial area} $a(F)$ of $F$ to be the
sum of the combinatorial areas of its 0-handles.

The reason for the above terminology is an analogy from
hyperbolic geometry. In that context, the area
of a planar convex hyperbolic polygon is the sum of its
exterior angles, with $2\pi$ then subtracted.
In the above formula, there is an extra term
involving the number of intersections between
$N$ and $\partial M$. Again, this can be understood
intuitively from hyperbolic geometry. If $M - \partial M$ is
hyperbolic, then the parts of $F$ approaching a
cusp of $M - \partial M$ most naturally inherit zero interior angle
there.

Note that there is only a finite number of normal curves
in any given 0-handle of ${\cal H}$ up to ambient isotopy
leaving ${\cal H}$ invariant. Thus,
it is a simple process to check whether a given assignment
of angles actually gives rise to an angled spine.
Note also that the condition that a particular choice of
exterior angles gives rise to an angled spine is a `convex
condition'. In other words, given two choices of angles,
each of which yields an angled spine, then any convex linear combination
of these angles also gives an angled spine. The following
is a particularly important example of an angled spine.

\noindent {\bf Lemma 4.1.} {\sl An angled ideal
triangulation determines an angled spine.}

\noindent {\sl Proof.} 
The interior angles around each edge of the ideal triangulation
sum to $2 \pi$, which gives the second condition in the definition
of an angled spine.
To see the first condition, observe that each component of ${\cal F}$
is as shown in Figure 10. Each normal curve $N$
which intersects $\partial M$ in more than one arc
automatically has non-negative combinatorial area. The normal curves
intersecting $\partial M$ in at most one arc are given
in Figure 13 (up to the obvious symmetries of ${\cal F}$), 
and can be seen to have non-negative combinatorial area. 
$\square$

\vskip 12pt
\centerline{\psfig{figure=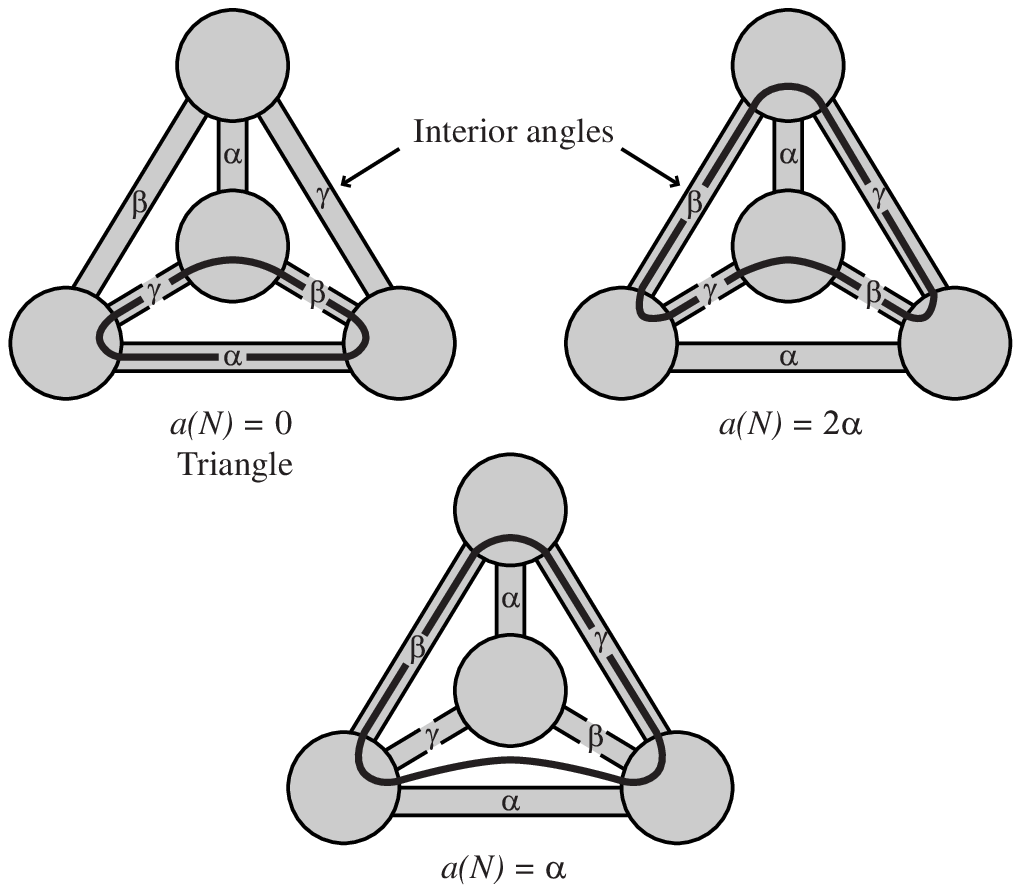}}
\vskip 12pt
\centerline{Figure 13.}

The above proof has the following corollary. Let
${\cal H}$ be a handle structure dual to an angled ideal
triangulation. Then the only normal curves with zero combinatorial area in
a 0-handle of ${\cal H}$ are {\sl triangles} (as in the leftmost
diagram of Figure 13) and {\sl boundary bigons},
one of which is shown in Figure 14. 

\vskip 12pt
\centerline{\psfig{figure=bbigon.ps}}
\nobreak
\vskip 6pt
\nobreak
\centerline{Figure 14.}

\noindent More generally, we will say that
a normal curve $N$ in the boundary of a 0-handle of an
angled spine (not necessarily arising from an angled ideal
triangulation) is a {\sl boundary bigon} if it
is disjoint from ${\cal F}^1$ and encloses a disc in
$\partial {\cal H}^0$ containing a single 1-handle of ${\cal F}$.

We wish to consider surfaces $F$ in $M$ which need not be
embedded; in this case, the theory of normal surfaces is less well
understood. In fact, we will consider surfaces $F$ which may even have
boundary components not lying in $\partial M$.
But, we will insist that $\partial F$ is
disjoint from ${\cal H}^2$, that its intersection
with ${\cal H}^1$ respects the product structure on
${\cal H}^1$, and that if a component of $\partial F$ touches
$\partial M$ then it lies entirely in $\partial M$. 
It is easy to see that any such surface $F$
can be ambient isotoped (keeping $\partial F$ invariant) 
so that, afterwards, its
intersection with any $i$-handle $D^i \times D^{3-i}$ is
of the form $D^i \times G$, for some codimension one 
submanifold $G$ of $D^{3-i}$. After a further
homotopy (keeping $\partial F$ fixed) we may assume
that $F$ intersects each handle in a collection of discs,
providing that $F$ is homotopically 
incompressible and has no 2-sphere components, and $M$ is irreducible.
Each 0-handle $H$ of $F$ intersects $\partial {\cal H}^0$
in a 1-manifold. This 1-manifold is a collection
of arcs if $H$ touches $\partial F - \partial M$,
but otherwise it is a closed curve. It need
not be embedded, but we can ensure (after a further
homotopy keeping $\partial F - \partial M$ fixed) that each arc
of intersection with ${\cal F}^0$, ${\cal F}^1$ and
$\partial {\cal H}^0 - {\cal F}$ is embedded. If $F$ is
homotopically $\partial$-incompressible and no component
is a disc parallel to a disc in $\partial M$, and $\partial M$
is incompressible, then we may further homotope $F$ 
(keeping $\partial F - \partial M$ fixed)
so that each curve of $F \cap \partial {\cal H}^0$ satisfies
conditions (i) - (iv) in the above definition of normality
(with the possible exception of arc components of 
$F \cap \partial {\cal H}^0$ which are allowed to lie entirely
in ${\cal F}^0$).
We say that $F$ is {\sl admissible} if it satisfies each
of the above restrictions. Also, we say that a curve $N$
in $\partial {\cal H}^0$ is {\sl admissible} if it satisfies
conditions (i) - (iv) and each arc of intersection with
${\cal F}^0$, ${\cal F}^1$ and $\partial {\cal H}^0 - {\cal F}$ is
embedded. We will not insist that (v) and (vi) hold. 

We wish to define the combinatorial area of an admissible
surface $F$. Let $H$ be a 0-handle of $F$, and, as above, 
let $\epsilon_1, \dots, \epsilon_n$
be the exterior angles of the 1-handles of ${\cal F}$ over which
$\partial H$ runs (counted with multiplicity if it runs over a 1-handle
of ${\cal F}$  more than once). Then define the {\sl combinatorial area}
of $H$ to be
$$a(H) = (\sum_{i=1}^n \epsilon_i) -2 \pi +
\pi \vert H \cap \partial M \vert +
3 \pi \vert \partial H - \partial {\cal H}^0 \vert.$$
The coefficient $3 \pi$ is somewhat arbitrary; it was chosen
so that the combinatorial area of $H$ is automatically positive
if $\partial H \cap \partial {\cal H}^0$ is not a closed curve.
Define the {\sl combinatorial area} of $F$ to be the sum of the combinatorial areas
of its 0-handles. For combinatorial area to have useful applications to
non-embedded surfaces, we need the following fact.

\noindent {\bf Lemma 4.2.} {\sl Let ${\cal H}$ be a thickened
angled spine. Then each closed admissible curve $N$
in $\partial {\cal H}^0$ has non-negative combinatorial area.
Furthermore, if $N$ is not a normal curve,
then its combinatorial area is positive.}

\noindent {\sl Proof.} We will prove this by induction on
the number of arcs of intersection between $N$ and ${\cal F}^1$.
If $N$ is normal then, from the definition of an angled spine, its
combinatorial area is non-negative. Suppose therefore
that $N$ is not normal and that its combinatorial area is non-positive.

The induction starts with $N$ disjoint from ${\cal F}^1$.
Since $N$ is not normal, it must have at least two arcs of intersection with
$\partial M$. This gives it non-negative area. Since
we are assuming that the area of $N$ is non-positive, it must be
zero. Hence, $N$ has precisely two arcs of intersection
with $\partial M$. If $N$ is non-embedded in some
0-handle of ${\cal F}$, then it is not hard to construct
a normal curve which intersects $\partial M$
in a single arc and which is disjoint from ${\cal F}^1$.
But this normal curve has negative combinatorial area,
contradicting the assumption that this is an angled spine.
Hence, the only way that $N$ can fail to be normal
if it is non-embedded in $\partial M$ or violates condition (vi) of 
the definition of normality. Again, we can construct
a negative area normal curve in $\partial {\cal H}^0$.
This starts the induction.  

Now suppose that $N$ intersects ${\cal F}^1$.
If $N$ violates condition (vi), then it must have at least
two arcs of intersection with $\partial M$. This gives it non-negative area. Also, 
if its area is zero, then these are the only arcs of intersection with 
$\partial M$, and also $N$ is disjoint from ${\cal F}^1$. 
But we have already dealt with this case.

If $N$ violates condition (v) of the
definition of normality, it contains two sub-arcs 
in the same 1-handle of ${\cal F}$. 
Extend these to maximal sub-arcs $A_1$ and $A_2$ of $N$,
with the property that each $A_i$ lies in ${\cal F}$ and
that $A_1$ and $A_2$ run along the same handles of ${\cal F}$.
Since we are assuming that $N$ has non-positive
area and that it runs over some 1-handle of ${\cal F}$,
it therefore intersects $\partial M$ in
at most one arc. Hence, at least one arc $A_3$ of
${\rm cl}(N - (A_1 \cup A_2))$ lies in ${\cal F}$.
If $\partial A_3$ lies at the same ends of
$A_1$ and $A_2$, then $A_3$ can be closed up to
form an admissible curve. If $\partial A_3$ lies
at opposite ends of $A_1$ and $A_2$, then the
$A_1 \cup A_3$ can be closed to form an admissible curve.
In either case, the inductive
hypothesis gives that the combinatorial area of this
admissible curve is non-negative and hence the combinatorial area
of $N$ is positive. 

If $N$ is not embedded but does not violate (v) or (vi), 
then let $P_1$ and $P_2$ be points on $N$
which are coincident in $\partial {\cal H}^0$. Then $P_1$ and
$P_2$ lie in some 0-handle of ${\cal F}$.  If both arcs of
$N - (P_1 \cup P_2)$ meet $\partial M$, then, since $N$
has non-positive area, it is disjoint from ${\cal F}^1$. We have proved
the lemma in this case. Suppose therefore that at least one
arc of $N - (P_1 \cup P_2)$ lies entirely in ${\cal F}$.
Join the endpoints of this arc to form a closed curve $N_1$.
Then $N_1$ satisfies all the conditions of normality,
except that it may be non-embedded. However, $N_1$ has
fewer self-intersections than $N$, and so we may repeat
this process and end with a normal curve. This has
non-negative area and so $N$ has positive area. $\square$

The following result is analogous to the
Gauss-Bonnet theorem.

\noindent {\bf Proposition 4.3.} {\sl Let $F$ be an
admissible surface in a handle structure ${\cal H}$
arising from an angled spine. Let ${\rm Length}(\partial F - \partial M)$ 
be the number of arcs of intersection between $\partial F - \partial M$
and the 1-handles of ${\cal H}$. Then the combinatorial area 
of $F$ satisfies
$a(F) = -2\pi\chi(F) + 2 \pi {\rm Length}(\partial F - \partial M)$.}

\noindent {\sl Proof.} This is a very straightforward
counting argument. Examine the formula for combinatorial
area, term by term. Summing the $(-2 \pi)$ terms for
each 0-handle of $F$ gives $-2 \pi \vert {\cal H}^0 \cap F \vert$.
Each exterior angle is $\pi$ minus the corresponding interior
angle. The interior angles sum to $2 \pi \vert {\cal H}^2 \cap F \vert$.
We need to know the total number $\vert F \cap {\cal F}^1 \vert$ of exterior angles
of $F$. Arranged around each 0-handle $H$ of $F$ are
components of $\partial H \cap {\cal F}^0$, alternating
with components of $\partial H \cap {\cal F}^1$, $\partial H \cap
\partial M$
and $\partial H - \partial {\cal H}^0$. Each component of
$F \cap {\cal H}^1$ gives two components of $F \cap {\cal F}^0$.
Hence,
$$2 \vert {\cal H}^1 \cap F \vert = \vert F \cap {\cal F}^0 \vert
= \vert F \cap {\cal F}^1 \vert + \vert F \cap \partial {\cal H}^0 \cap \partial M \vert
+ {\rm Length}(\partial F - \partial M).$$
Therefore, $\pi \vert F \cap {\cal F}^1 \vert$, together with the sum of the
$3 \pi \vert \partial H - \partial {\cal H}^0 \vert$
and $\pi \vert H \cap \partial M \vert$ terms gives
$2 \pi \vert {\cal H}^1 \cap F \vert + 2\pi {\rm Length}(\partial F
- \partial M)$. $\square$

\vfill\eject
Casson used this to show the following result.

\noindent {\bf Proposition 4.4.} {\sl An orientable 3-manifold $M$
with an angled ideal triangulation contains no normal 2-spheres and the
only normal tori arise as links of the ideal vertices.
Also, each boundary component of $M$ is a torus.}

\noindent {\sl Proof.} Let $F$ be a normal 2-sphere or torus.
Then the combinatorial area of each 0-handle of $F$ is non-negative,
and so the combinatorial area of $F$ is non-negative.
But, by Proposition 4.3, the combinatorial area of a 2-sphere is
$-4 \pi$ and that of a torus is zero. Hence, we can have no
normal 2-sphere. Each 0-handle of a normal torus has
zero combinatorial area, and therefore, by the remarks after Lemma 4.1,
is a triangle as in Figure 13. These
triangles join to form the link of an ideal vertex. 
Furthermore, the link of an ideal vertex is comprised
of triangles, each with zero combinatorial area. Hence,
by Proposition 4.3, each boundary component of $M$ has
zero Euler characteristic, and so is a torus. $\square$

In a similar spirit, we can restrict the possible
surfaces with non-negative Euler characteristic in
an angled spine.

\noindent {\bf Proposition 4.5.} {\sl Let $M$ be a
3-manifold with a handle structure ${\cal H}$ arising from an angled
spine. Then ${\cal H}$ contains no normal properly embedded spheres or
discs. Also, if the angled spine is dual to an angled ideal
triangulation, then each normal annulus properly embedded in $M$ is constructed by
picking a 2-handle $H_2$ of ${\cal H}$,
and then taking ${\rm cl}(\partial {\cal N}(H_2) - \partial M)$,
where ${\cal N}(H_2)$ is a small regular neighbourhood of
$H_2$.}

\noindent {\sl Proof.} As explained above, the combinatorial
area of a normal surface in non-negative. However,
by Proposition 4.3, the combinatorial area of a properly
embedded disc or sphere is negative. Thus, ${\cal H}$ can contain
no normal spheres or properly embedded discs. Similarly, the combinatorial area
of a properly embedded annulus $F$ is zero. So, if ${\cal H}$ is
dual to an angled ideal triangulation, each 0-handle
of $F$ is either a triangle or a boundary bigon.
At least one 0-handle of $F$ is a boundary bigon.
Also, a boundary bigon and a triangle cannot be
adjacent 0-handles in $F$, and therefore, every 0-handle
of $F$ is a boundary bigon. These combine to form an
annulus as described in the proposition. $\square$

\vfill\eject
\noindent {\bf Corollary 4.6.} {\sl An orientable 3-manifold with an
angled ideal triangulation is irreducible, atoroidal and not
Seifert fibred.}

\noindent {\sl Proof.} Irreducibility and atoroidality
follow from Proposition 4.4. Seifert
fibred manifolds with non-empty boundary are either the
solid torus or have essential properly embedded annuli.
In the former case, there is a meridian disc in normal form.
In the latter case, we may find an essential properly embedded
annulus in normal form. This contradicts Proposition 4.5. $\square$

We cannot in general deduce information about annuli
and tori in angled spines, since there may exist
many types of normal curves in $\partial {\cal H}^0$
with zero combinatorial area.
The following result asserts that manifolds with angled spines are
fairly ubiquitous.

\noindent {\bf Theorem 4.7.} {\sl Let $M$ be a compact orientable
irreducible atoroidal 3-manifold with non-empty boundary. Suppose that
$M$ contains no properly embedded essential discs or annuli, 
and that $M$ is not a 3-ball. Then $M$ has an angled spine.}

\noindent {\sl Proof.} Let us first deal with the case where each
boundary component of $M$ is a torus. Then, by [24], $M - \partial M$ 
has a complete finite volume hyperbolic structure. By [7], $M - \partial M$
admits an expression as a union of convex hyperbolic ideal polyhedra with faces
identified isometrically in pairs. The dual of this is a spine, with
an interior angle associated with each component of ${\cal F}^1$,
arising from the interior dihedral angle of the
corresponding edge in one of the polyhedra.
We claim that this is actually an angled spine. As in the proof
of Lemma 4.1, the second part of the definition of an angled spine is
immediately verified. To check the first part, consider a normal
curve $N$ in $\partial {\cal H}^0$. We need to show that $a(N) \geq 0$.
If $N$ lies entirely in ${\cal F}$, then a result of Rivin (Theorem 1 
of [21]) about convex
hyperbolic ideal polyhedra implies that the sum of the exterior angles
of $N$ is at least $2 \pi$, and so $a(N) \geq 0$. However,
$N$ may intersect $\partial M$.
If it intersects $\partial M$ in at least two arcs, then
$a(N) \geq 0$. Suppose that $N$ intersects $\partial M$ in precisely
one arc $A$. If we replace $A$ with an arc $A'$ in ${\cal F}$
skirting around the component of $\partial {\cal H}^0 - {\cal F}$
containing $A$, the result is an admissible curve $N'$ in ${\cal F}$. There
are two choices $A'_1$ and $A'_2$ for $A'$ corresponding to the two ways of going
around the component of $\partial {\cal H}^0 - {\cal F}$.
The exterior angles of $A'_1$ and $A'_2$ sum to $2 \pi$,
and hence we can ensure that the exterior angles of $A'$
sum to at most $\pi$. Hence, $a(N) \geq a(N')$. Now, the curve
$N'$ might not be normal, but by the argument of Lemma 4.2,
$a(N')$ is at least that of a normal curve in ${\cal F}$, and
so $a(N') \geq 0$. This verifies that we have an angled spine. 

Consider now the case where $M$ contains some non-toral boundary
components. Let $T$ be the toral boundary components of $M$ (possibly,
$T = \emptyset$). Let $Y$ be the 3-manifold obtained by doubling $M - T$
along $\partial M - T$. Then our assumptions about $M$ imply that
$Y$ is irreducible and atoroidal. Now, $Y$ is Haken, since $\partial M -
T$ is incompressible. Hence, $Y$ is not a Seifert fibre space, since
otherwise, according to VI.34 of [14], $\partial M - T$ would separate
$Y$ into two $I$-bundles over compact surfaces, which would imply
that $M$ contained a properly embedded essential annulus, contrary to
assumption. Hence, by [24], $Y$ has a complete finite volume hyperbolic structure. It admits
an isometric involution $\tau$ which swaps the two copies of $M - T$.
The fixed point set of $\tau$ is a totally geodesic copy of $\partial M
- T$. Hence, we have the well-known result that $M$ admits
a complete hyperbolic structure with $\partial M - T$ being 
totally geodesic.

In the case where $T \not= \emptyset$, $Y$ admits a canonical expression as a 
union of convex hyperbolic ideal polyhedra with faces identified in pairs [7]. 
This is preserved by $\tau$. Hence, the intersection of these polyhedra with
one half of $Y$ then gives an expression of $M - T$ as a union of convex
hyperbolic polyhedra $P$ (with some vertices possibly at infinity), with
$\partial M - T$ being a union of faces of $P$. Each edge
touching $\partial M - T$ either lies entirely in $\partial M - T$
or is at right-angles to it. 

Similarly, when $T = \emptyset$,
according to [16], $M$ admits such a representation as a union of convex hyperbolic
polyhedra $P$, with each edge touching $\partial M - T$ either lying entirely 
in $\partial M - T$ or at right-angles to it. These polyhedra are
known as `truncated polyhedra' since they come from hyper-ideal
polyhedra by performing a perpendicular truncation of the
infinite-volume ends.

The angled spine of $M$ is obtained by dualising $P$ 
in the following fashion: associate a 0-handle of the spine
with each polyhedron of $P$, associate a 1-handle of the
spine with each face of $P - \partial M$, and associate a 2-handle of
the spine with each edge of $P$ not lying entirely in $\partial M$.
Again, to show that this is an angled spine, we must show that each
normal curve $N$ in $\partial {\cal H}^0$ has non-negative combinatorial
area. As in the case above where $P$ is ideal, we can assume that 
$N$ misses $T$ (but it might run over components of $\partial M - T$).
Hence, it corresponds to a normal curve (also called $N$)
in boundary of a polyhedron $P'$ of $P$ which is
transverse to the edges of $P'$.
Each arc of intersection with $\partial M - T$ contributes
two external angles of $\pi/2$ to $N$. Hence, we may assume that $N$ has
at most one arc of intersection with $\partial M - T$.
If it misses $\partial M - T$, then it corresponds to
a path $N_1$ in the boundary of the hyper-ideal polyhedron $P_1$ obtained by
re-attaching an infinite volume end to each component
of $P' \cap (\partial M - T)$. If $N$ hits $\partial M - T$,
then let $P_2$ be the hyper-ideal polyhedron obtained
by doubling $P'$ along the component of $P' \cap (\partial M - T)$
intersecting $N$, and then attaching infinite volume ends
to the remaining components of $P' \cap (\partial M - T)$.
Two copies of $N$ (one in each half of $P_2$) join to form a closed 
curve $N_2$ in the boundary of $P_2$, with $a(N_2) = 2 a(N)$.
By the argument of [21] generalised to hyper-ideal polyhedra,
$N_1$ and $N_2$ have non-negative combinatorial area. Hence, so
does $N$. $\square$

The combinatorial $2 \pi$ theorem outlined in the
introduction will be proved using results about
admissible surfaces with boundary.
Let ${\cal H}$ be the handle structure arising from
an angled spine. We wish to
define the combinatorial length of a slope $s$ on $\partial M$.
Pick a curve $C$ representing a non-zero multiple of $s$ which respects the
handle structure on $\partial M$, in the sense that it is disjoint
from the 2-handles and is vertical in the 1-handles.
Let $C^0_1, \dots, C^0_n$ be
the sequence of arcs of intersection between
$C$ and the 0-handles of $\partial M$, and let
$C^1_1, \dots, C^1_n$ be the
sequence of arcs of intersection between $C$
and the 1-handles of $\partial M$, where $C^1_i$ lies between
$C^0_i$ and $C^0_{i+1}$.
We view $C^0_{n+i}$ as $C^0_i$, and $C^1_{n+i}$ as $C^1_i$.
Define an {\sl inward extension} of $C$ to be
a surface ${\cal E}$ lying in $M$, comprised of handles
$H^0_1, \dots, H^0_n, H^1_1, \dots, H^1_n$, such that
\item{$\bullet$} each $H^j_i$ is a disc in a $j$-handle of ${\cal H}$;
\item{$\bullet$} each $H^j_i$ contains $C^j_i$;
\item{$\bullet$} each $H^1_i$ respects the product structure
of the 1-handle of ${\cal H}$ in which it lies;
\item{$\bullet$} each intersection $H^0_i \cap \partial {\cal H}^0$ is admissible;
\item{$\bullet$} the arcs of ${\cal F}^0 \cap H^0_i$
and ${\cal F}^0 \cap H^1_i$ touching $C^0_i \cap C^1_i$ agree;
\item{$\bullet$} the arcs of ${\cal F}^0 \cap H^1_i$
and ${\cal F}^0 \cap H^0_{i+1}$ touching $C^1_i \cap C^0_{i+1}$ agree.

\noindent Define the {\sl weight} of
${\cal E}$ to be
$$w({\cal E}) = \sum_{i=1}^n a(H^0_i) / \vert H^0_i \cap \partial M \vert.$$
Define the {\sl combinatorial length} of $C$ to be
$$l(C) = \inf \lbrace w({\cal E}) : \hbox{${\cal E}$ is an inward extension of
$C$} \rbrace,$$
and define $l(C)$ as infinite if $C$ has no inward extension. This
happens if an arc $C_i^0$ has endpoints in the same component
of $\partial {\cal F}^0 - {\cal F}^1$, for then a curve $\partial H^0_i$
containing $C^0_i$ must fail condition (iv) in the definition of normality,
and so cannot be admissible. In fact, if no arc $C_i^0$ is
of this form, then $C$ has an inward extension which is an annulus
formed by pushing $C$ a little into the interior of $M$.
Let the {\sl combinatorial length} of $s$ be
$$l(s) = \inf \lbrace l(C) : \hbox{$C$ is a curve representing 
a non-zero multiple of $s$} \rbrace.$$

These definitions have been designed specifically so that the
following proposition holds.

\noindent {\bf Proposition 4.8.} {\sl Let $F$ be
an admissible surface in a handle structure on $M$
arising from an angled spine. Let $C_1, \dots, C_m$ be
the components of $\partial F \cap \partial M$, each $C_j$
representing a non-zero multiple of some slope $s_{i(j)}$. Then
$$a(F) \geq \sum_{j=1}^m l(s_{i(j)}).$$}

\noindent {\sl Proof.} The handles of $F$ touching $C_j$
form an inward extension of $C_j$. Summing the weights of
these inward extensions over all $C_j$ gives the required
inequality. $\square$

The following result is our `word hyperbolic Dehn surgery
theorem' for angled spines. Note that this, together with Corollary 4.6,
gives the version of Theorem 4.9 mentioned in the introduction.

\noindent {\bf Theorem 4.9.} {\sl Let $M$ be a compact orientable
3-manifold with an angled spine. Suppose that $M$ is atoroidal and not
a Seifert fibre space, and has boundary a non-empty union of tori.
Let $s_1, \dots, s_n$ be a collection of slopes on $\partial M$, 
with one $s_i$ on each component of $\partial M$, and each
with combinatorial length more than $2\pi$.
Then the manifold obtained by Dehn filling $M$
along these slopes is irreducible, atoroidal and not Seifert fibred,
and has infinite, word hyperbolic fundamental group.}

\noindent {\sl Proof.} Note that $M$ is irreducible, since
it has an angled spine. Since we are assuming that it
is atoroidal and not a Seifert fibre space, its interior
supports a complete finite volume hyperbolic structure [24].
Suppose first that $M(s_1, \dots, s_n)$
is reducible or toroidal. Then we can find $F$ a punctured sphere
or torus in $M$, with $F$ incompressible and $\partial$-incompressible
in $M$. Hence, $F$ can be ambient isotoped into normal form in
${\cal H}$, the thickened angled spine of $M$.
Let $C_1, \dots, C_{\vert \partial F \vert}$ be the boundary
components of $F$. We can ensure that $C_j$ is
essential in $\partial M$ having slope
$s_{i(j)} \in \lbrace s_1, \dots, s_n \rbrace$. Proposition 4.3 gives that
$$a(F) = -2\pi \chi(F) \leq 2\pi \vert \partial F \vert,$$
but Proposition 4.8 gives that
$$a(F) \geq \sum_{j=1}^{\vert \partial F \vert} l(s_{i(j)})
> 2 \pi \vert \partial F \vert,$$
which together give a contradiction. A similar
argument gives that the core of each surgery solid
torus in $M(s_1, \dots, s_n)$ has infinite order
in $\pi_1(M(s_1, \dots, s_n))$. 

To show the word hyperbolicity of
$\pi_1(M(s_1, \dots, s_n))$, we will use Gabai's
Ubiquity theorem. Let $K$ be a curve in $M$
which is homotopically trivial in $M(s_1, \dots, s_n)$.
Assume that $K$ is disjoint from the 2-handles of ${\cal H}$ and
respects the product structure on the 1-handles.
Then, there is a compact planar surface $F$ in $M$,
with $\partial F$ being $K$ together with
curves $C_1, \dots, C_{\vert F \cap \partial M \vert}$ on $\partial M$, each
$C_j$ being a non-zero multiple $k_j$ of some slope $s_{i(j)}$.
By Lemma 2.5, we may assume that $F$ is homotopically incompressible
and homotopically $\partial$-incompressible. We may therefore homotope
$F$ (keeping $K$ invariant) to an admissible surface.
Let ${\rm Length}(K)$ be the number of arcs of intersection between
$K$ and ${\cal H}^1$, and let
$\epsilon > 0$ be a real number with $l(s_i) \geq 2\pi+ \epsilon$
for all $i$. Then
Propositions 4.8 and 4.3 give that
$$\eqalign{(2 \pi + \epsilon) \vert F \cap \partial M \vert
&\leq \sum_{j=1}^{\vert F \cap \partial M \vert} l(s_{i(j)}) \cr
&\leq a(F) = -2\pi \chi(F) + 2\pi {\rm Length}(K) \cr
&< 2\pi \vert F \cap \partial M \vert + 2\pi {\rm Length}(K),\cr}$$
from which we get
$$\vert F \cap \partial M \vert < (2 \pi / \epsilon)
{\rm Length}(K).$$
Theorem 2.1 gives us that $\pi_1(M(s_1, \dots, s_n))$
is word hyperbolic. $\square$

The conclusions of Theorem 4.9 still hold (with the exception of
word hyperbolicity) if one leaves some boundary components of
$M$ unfilled. Also, one can allow $M$ to have some non-toral
boundary components.

Of course, the above theorem begs the question of which slopes
on $\partial M$ have combinatorial length more than $2\pi$.
The purpose of the following proposition is to give
a simple lower bound in the case of angled ideal triangulations.

\noindent {\bf Proposition 4.10.} {\sl Let $M$ be an orientable 3-manifold with
an angled ideal triangulation. Define the length of
each edge $E$ of the associated triangulation of $\partial M$
to be $\min \lbrace \alpha_1, \dots, \alpha_6 \rbrace /2 $,
where $\alpha_1, \dots, \alpha_6$ are the interior
angles of the two triangles adjacent to $E$. Each simplicial curve
in the 1-skeleton of $\partial M$ then inherits a length. The shortest length
of a curve of slope $s$ is a lower bound for the combinatorial length $l(s)$.}

\noindent {\sl Proof.} The ideal triangulation of $M$ dualises
to a special spine, and so the boundary of $M$ inherits a
handle structure in which each 0-handle has valence three.
Consider a curve $C$ representing a non-zero multiple of the slope $s$ 
in $\partial M$, such that $C$ is disjoint
from the 2-handles and respects the product structures
on the 1-handles. As above, let $C^1_1, \dots, C^1_n$ be the
arcs of intersection between $C$ and the 1-handles of $\partial M$.
For each $C^1_i$, let $E_i$ be the associated edge in
the dual triangulation. Let ${\cal E}$ be an inward
extension of $C$.

We term $C^1_i$ a {\sl hugging arc} if $C^1_{i-1}$, $C^1_i$ and
$C^1_{i+1}$ are all adjacent to some 2-handle of
$\partial M$ and are arranged in order around
that 2-handle, as shown in Figure 15.
The point behind this definition is that
if $H^0_i$ and $H^0_{i+1}$ are both
boundary bigons in the inward extension ${\cal E}$, then
$C^1_i$ must be a hugging arc.

\vskip 12pt
\centerline{\psfig{figure=hugarc.ps}}
\vskip 12pt
\centerline{Figure 15.}

\noindent {\sl Claim.} Suppose that $H = H^0_i$ is not a boundary
bigon. Let $\lbrace \alpha_1, \alpha_2, \alpha_3 \rbrace$ be the interior
angles of the 0-handle of ${\cal H}$ containing $H$.
Then $a(H) / \vert H \cap \partial M \vert \geq
\min \lbrace \alpha_1, \alpha_2, \alpha_3 \rbrace$.

Let $\epsilon_1, \dots, \epsilon_t$
be the exterior angles inherited by $\partial H$.
If $\partial H$ intersects $\partial M$
in at least three arcs, then
$$\eqalign{
{a(H) \over \vert H \cap \partial M \vert}
&= {(\sum_{j=1}^t \epsilon_j) -2 \pi +
\pi \vert H \cap \partial M \vert + 3\pi \vert \partial H - \partial
{\cal H}^0 \vert
\over \vert H \cap \partial M \vert} \cr
&\geq \pi \Big( 1 - {2 \over \vert H \cap \partial M \vert} \Big)
\geq \pi/3.}$$
But $\alpha_1 + \alpha_2 +\alpha_3 = \pi$, and therefore
$\min \lbrace \alpha_1, \alpha_2, \alpha_3 \rbrace \leq \pi/3.$
Therefore, in the case where $\vert H \cap \partial M \vert \geq 3$,
the claim is proved. For $\vert H \cap \partial M \vert = 2$, we have
$$\eqalign{
{a(H) \over \vert H \cap \partial M \vert} &= {(\sum_{j=1}^t \epsilon_i) -2 \pi +
\pi \vert H \cap \partial M \vert + 3 \pi \vert \partial H - \partial {\cal H}^0
\vert \over \vert H \cap \partial M \vert} \cr
&\geq \min \lbrace \pi - \alpha_1, \pi - \alpha_2, \pi - \alpha_3, 3\pi \rbrace
/2 \cr
&\geq \min \lbrace \alpha_1, \alpha_2, \alpha_3 \rbrace.\cr}$$
For $\vert H \cap \partial M \vert = 1$ and $\partial H$ not normal,
then the proof of Lemma 4.2 gives that the combinatorial area of $H$ is at least
$\pi$. The normal curve $N$ with $\vert N \cap \partial M \vert = 1$
is given in Figure 13, and its combinatorial area is 
at least $\min \lbrace \alpha_1, \alpha_2, \alpha_3 \rbrace$. 
This proves the claim.

We will show that $C$ naturally determines a
curve $C'$ in the 1-skeleton of the dual triangulation, in the following manner.
Keeping the midpoints of $C^1_i$ and $C^1_{i+1}$ fixed,
homotope the arc of $C$ lying between them, so that
it runs along $E_i$ and then back up $E_{i+1}$.
(We can do this because each 0-handle of the
handle structure on $\partial M$ has valence three.)
The curve which we have constructed
lies in the 1-skeleton of the triangulation
of $\partial M$, but need not be simplicial, since the vertex joining
$E_{i-1}$ to $E_i$ may be the same as the vertex
joining $E_i$ to $E_{i+1}$. This happens precisely when
$C^1_i$ is a hugging arc. However, a further
homotopy creates a simplicial curve $C'$,
with each edge $E_i$ of $C'$ being dual to a
non-hugging arc $C^1_i$.
Let $\lbrace \alpha_1, \dots, \alpha_6 \rbrace$ be the six
interior angles of the two triangles adjacent to $E_i$.
Then, at least one of $\partial H^0_i$ and $\partial H^0_{i+1}$
is not a boundary bigon. Therefore, by the claim,
$$
\left( {a(H^0_i) \over \vert H^0_i \cap \partial M \vert}
+ {a(H^0_{i+1}) \over \vert H^0_{i+1} \cap \partial M \vert} \right) /2 \geq
\min \lbrace \alpha_1, \dots, \alpha_6 \rbrace/2.$$
Summing this from $i=1$ to $n$ gives that the combinatorial length of $C$ is at
least the length of $C'$ in the path metric on the 1-skeleton
of the triangulation of $\partial M$. Since $C'$ represents
a non-zero multiple of the slope $s$, it is possible to perform
a cut-and-paste on $C'$ in the vertices of $\partial M$, making it 
into a collection of curves, each with slope $s$. In particular,
the length of $C'$ is at least the length of the shortest curve
in the 1-skeleton of $\partial M$ with slope $s$. $\square$

Note that the factor $1/2$ arises in the definition of edge
length in the triangulation of $\partial M$ to avoid
double counting. This is because we may have successive
edges $E_i$ and $E_{i+1}$ in $C'$, dual to non-hugging
arcs $C^1_i$ and $C^1_{i+1}$, but with $\partial H^0_i$ and
$\partial H^0_{i+2}$ both boundary bigons. Therefore,
we must share the combinatorial area of $H^0_{i+1}$ out
between $E_i$ and $E_{i+1}$. However,
if successive edges $E$ and $E'$ of $C'$ are not dual to 
successive edges of $C$ (because there is some hugging
arc between them), then this double counting cannot occur.
Hence, the lower bound in Proposition 4.10
need not be sharp. Nevertheless, Proposition 4.10 has
the following immediate corollary.

\noindent {\bf Corollary 4.11.} {\sl Let $M$ be an orientable
3-manifold with an angled ideal triangulation. Then
at most finitely many slopes on each component of $\partial M$
have combinatorial length no more than $2 \pi$.}

\vfill\eject
\centerline{\caps 5. Surgery along alternating links}
\vskip 6pt

Our goal here is to use the techniques of the previous section
to deduce the following word hyperbolic Dehn surgery theorem for 
alternating links. 

\noindent {\bf Theorem 5.1.} {\sl Let $D$ be a connected
prime alternating diagram of a link $L$ in $S^3$. For each
component $K$ of $L$, pick a surgery coefficient $p/q$ (in its lowest terms)
with $\vert q \vert > 8/t(K,D)$. The manifold
obtained by Dehn surgery along $L$, via these surgery
coefficients, is irreducible, atoroidal and not a Seifert fibre space, and has
infinite, word hyperbolic fundamental group.}

The number $t(K,D)$ was defined in the introduction to
be half the number of non-bigon edges of $G(D)$ lying in $K$.
Similarly, the twist number $t(D)$ was defined to be half the number
of non-bigon edges in $G(D)$. It has the following simple interpretation.
Define a {\sl twist} of $D$ to be a maximal connected subgraph of $G(D)$
with edges comprised of bigon edges. A single vertex of $G(D)$ incident to
no bigon edge is viewed as a twist. Twists in a connected diagram $D$
can be of several forms: either that shown in Figure 16; or they can 
comprise the whole diagram, in which case $D$ is the standard
diagram of the $(2,n)$-torus link for some $n$; or they can form
a Hopf link summand with two crossings. Henceforth, suppose 
that $D$ has positive twist number and has no Hopf link summands 
with two crossings. Then, it is not hard to see
that $t(D)$ is simply the number of twists in $D$. 

\vskip 12pt
\centerline{\psfig{figure=atwist.ps}}
\vskip 12pt
\centerline{Figure 16.}

If one replaces each twist of $G(D)$ with a single
vertex, then one obtains a new 4-valent graph $G_1$ with precisely
$t(D)$ vertices. One can then reconstruct $G(D)$ from $G_1$
by replacing some or all of its vertices with twists.
This procedure can be effected in a surgical manner,
as follows. Start with $G_1$ and replace some or all of its vertices
with two vertices joined by a bigon region, creating a new 4-valent graph
$G_2$. Then assign under-over crossing information to $G_2$ creating
an alternating link $L_2$, say. Add unknotted components to
$L_2$, one component for each vertex of $G_1$, which encircles
the crossing (or pair of crossings) as shown in Figure 17.
Let $L_3$ be the resulting link. It is an augmented alternating 
link, in the sense of [1]. Then the link $L$ is obtained from
$L_3$, by performing $1/q$ surgery on each of the new unknotted
components (where $q \in {\Bbb Z}$ depends on the component of $L_3$).

\vskip 12pt
\centerline{\psfig{figure=augalt.ps}}
\vskip 12pt
\centerline{Figure 17.}

For a fixed integer $n$, there is only a finite number of
4-valent planar graphs $G_1$ with $n$ vertices. At each stage
of the construction of $L_3$ from $G_1$, we had only a finite number
of choices. Thus, we have shown that the alternating links
with a given twist number $n$ are all obtained from a finite
list of augmented alternating links via surgery on some of
the components of the link. Adams [1] showed that the complement of an augmented 
alternating link $L_3$ is hyperbolic. By [23], if the complement of the
alternating link $L$ is hyperbolic, then its
volume is at most that of the complement of $L_3$. Applying this
argument to knots with twist number at most 8, we have the
following corollary to Theorem 5.1.

\noindent {\bf Corollary 5.2.} {\sl There is a real number $c$ with
the property that if $K$ is an alternating knot whose complement has a
complete hyperbolic structure with volume at least $c$, 
then every non-trivial surgery on $K$ yields a
3-manifold which is irreducible, atoroidal and not Seifert fibred, and has
infinite, word hyperbolic fundamental group.}

It has been known for several years that the complement of a non-split
alternating link can be expressed as the union of two ideal
polyhedra glued along their faces. Dually, a non-split alternating
link complement has a spine with some nice properties. We
briefly recall here the salient points of the theory. A
more complete description can be found in [3] and [17].

The alternating link diagram $D$ is viewed as lying on a
2-sphere $S^2$ embedded in $S^3$. The link itself is
viewed as lying in $S^2$, except near each crossing,
where two arcs skirt above and below $S^2$. The spine
of the alternating link exterior has two 0-cells, one lying in each of the
two 3-balls $S^3 - {\rm int}({\cal N}(S^2))$; it has
a 1-cell for each region of the diagram, running between the two 0-cells; and it has
a 4-valent 2-cell for each crossing. Let ${\cal H}$ be the handle
structure arising by thickening the spine. The surface 
${\cal F} = {\cal H}^0 \cap ({\cal H}^1 \cup {\cal H}^2)$
has two components, one in each 0-handle of ${\cal H}$.
It is a crucial observation that each component
of ${\cal F}$ is dual to the underlying graph $G(D)$
of the diagram. Thus, dually, the link complement is obtained by
gluing together two ideal polyhedra $P_1$ and $P_2$, each with boundary
graph a copy of $G(D)$. 

\vskip 18pt
\centerline{\psfig{figure=altgraph.ps}}
\vskip 12pt
\centerline{Figure 18.}

Each face of $\partial P_1$ is
glued to the corresponding face of $\partial P_2$, but with
a twist. The orientation of the twist is determined
as follows. The complimentary regions of $G(D)$ admit
a black-and-white `checkerboard' colouring, with adjacent regions being
coloured with different colours. Then, we attach the black
(respectively, white) regions of $\partial P_1$ to the
corresponding region of $\partial P_2$ with a single
clockwise (respectively, anti-clockwise) twist. This
was first observed by Thurston [23] in the case of the
Borromean rings, where he likened this gluing procedure
to the gears of a machine. An example is given in Figure 19.

\vskip 12pt
\centerline{\psfig{figure=gears.ps}}
\vskip 12pt
\centerline{Figure 19.}

When the exterior $M$ of a non-split alternating link $L$
is given this handle structure, its boundary inherits a
particularly simple handle structure. The simplest way to see
this is to embed $\partial M$ as a normal surface. We then see that 
each component $T$ of $\partial M$ has the form of Figure 20.
Each 0-handle of $T$ has valence four. These 0-handles are arranged
along $T$, alternating between 0-handles above $S^2$ and
0-handles below $S^2$. If $C$ is a meridian curve in $\partial M$
encircling an edge $E$ of $G(D)$, then $C$ intersects two 1-handles
and two 0-handles of $\partial M$. We say that these handles are
{\sl associated} with the edge $E$. The two 1-handles associated with
$E$ each join the 0-handles associated with $E$. Hence, if $K$
is some component of $L$, the
total number of 0-handles in $\partial {\cal N}(K)$
is equal to the number of edges $e(K,D)$ of $K$ in $G(D)$. 

\vskip 12pt
\centerline{\psfig{figure=altcusp.ps}}
\vskip 12pt
\centerline{Figure 20.}

In order to utilise the results of Section 4, we
need to place an angled structure on the spine. 
The simplest way to do this is to assign an `exterior angle' in
the range $(0,\pi)$ to each edge of the graph $G(D)$. Since 
each 1-handle of ${\cal F}$ is dual to some edge of $G(D)$, 
this assigns an exterior angle to that 1-handle. Of course, one 
must check that these angles do
genuinely determine an angled spine. 

Recall that a diagram $D$ is {\sl reduced} if every crossing of $D$
is adjacent to four distinct complimentary regions of
$G(D)$. Note that if a diagram is prime and has more than
one crossing, then it is reduced. 

\noindent{\bf Proposition 5.3.} {\sl Assigning an exterior angle
to each edge of $G(D)$ gives an angled spine
if and only if each of the following conditions are satisfied:
\item{(i)} $D$ is reduced;
\item{(ii)} the exterior angles of the edges
around any vertex of $G(D)$ sum to $2 \pi$;
\item{(iii)} any simple closed curve in $S^2$ which
intersects $G(D)$ transversely in the interior of edges
of $G(D)$ and touches each edge at most
once, has total exterior angle at least $2 \pi$.

}

\noindent {\sl Proof.} The verification of the first part of the
definition of an angled spine is a simple argument analogous to the
first half of the proof of Theorem 4.7, and is omitted. To verify the
second part, recall that a 2-handle of ${\cal H}$ is
associated with each crossing of $D$. The four exterior angles
assigned to the edges of $G(D)$ emanating from that crossing
are precisely the four exterior angles assigned to the
2-handle. Hence, by (ii), the four exterior angles around that 2-handle sum
to $2 \pi$, and therefore so do the four interior angles. $\square$

In Figure 21, we assign exterior angles to the standard diagrams
of the figure-eight knot, the Borromean rings and the Whitehead link.
The angles satisfy the conditions of Proposition 5.3 and
so define angled spines on the link exteriors. In this case,
our choice of angles has been inspired by the hyperbolic structures
on these link complements.

\vskip 12pt
\centerline{\psfig{figure=angdiag.ps}}
\vskip 12pt
\centerline{Figure 21.}

It is the example of the Borromean rings which readily generalises
to {\sl any} connected prime alternating diagram $D$ with
more than one crossing.
If we assign an angle $\pi/2$ to each edge of $G(D)$, then the
conditions of Proposition 5.3 are satisfied. We will call this
the {\sl canonical angled spine} arising from $D$.

The major technical result of this section is the following.

\noindent {\bf Theorem 5.4.} {\sl Let $L$ be a link with a
connected prime alternating diagram $D$ having more than one
crossing. Give the exterior
of $L$ the canonical angled spine arising from $D$. Then
the combinatorial length of the slope $p/q$ on a
component $K$ of $L$ is at least $\vert q \vert \, t(K,D) \, \pi/4$.}

\noindent {\sl Proof of Theorem 5.1 from Theorem 5.4.}
If $D$ has precisely one crossing, then its twist number is
zero, and the theorem is trivial. So suppose that $D$
has more than one crossing.
Menasco [18] proved that the exterior of $L$ is irreducible, atoroidal
and not a Seifert fibre space,
unless $D$ is the standard diagram of a $(2,n)$-torus link for some $n$.
In this case, $t(K,D) = 0$, and the statement of the theorem
is empty. Otherwise, Theorem 5.1 follows from Theorem 5.4
and Theorem 4.9. $\square$

Let $M$ be the exterior of $L$.
Let $C$ be a curve on $\partial {\cal N}(K)$ representing a non-zero
multiple $k$ of the slope $p/q$. As before, we let $C^0_1, \dots, C^0_n$ be
the sequence of arcs of intersection between
$C$ and the 0-handles of $\partial M$, and let
$C^1_1, \dots, C^1_n$ be the
sequence of arcs of intersection between $C$
and the 1-handles of $\partial M$, where $C^1_i$ lies between
$C^0_i$ and $C^0_{i+1}$. Let ${\cal E} =
H_1^0 \cup \dots \cup H_n^0 \cup H_1^1 \cup \dots \cup H_n^1$ be an inward
extension of $C$.

Let $Y \rightarrow \partial {\cal N}(K)$ be the cover
corresponding to the subgroup of $\pi_1(\partial {\cal N}(K))$
generated by the meridian. In other words, unwind
$\partial {\cal N}(K)$ in the longitudinal direction.
We assign an $x$ co-ordinate to each 0-handle of $Y$,
with adjacent 0-handles having co-ordinates differing
by 1. Hence, $C$ lifts to an arc $\tilde C$ in $Y$ which starts at $x=0$
and ends at $x= k \vert q \vert e(K,D)$. Let $\tilde C^j_i$ be the 
lift of $C^j_i$ in $\tilde C$.
Let $x_i$ be the $x$ co-ordinate of $\tilde C^0_i$,
and let $x'_i = (x_{i-1} + x_{i+1})/2$.
Note that, from Figure 22, each $x'_i$ is an integer.

\vskip 12pt
\centerline{\psfig{figure=xcoord.ps}}
\vskip 12pt
\centerline{Figure 22.}

\noindent {\bf Lemma 5.5.} {\sl For each integer $j$ with 
$x'_1 \leq j < x'_1 + k \vert q \vert e(K,D)$, the 0-handles of $Y$ at
$x = j$ and $x = j+1$ have at least one arc $\tilde C^1_i$ running 
between them for which $x'_i \not= x'_{i+1}$.}

\noindent {\sl Proof.} Note that $\vert x'_{i+1} - x'_i \vert \leq 1$.
Also, if $x'_i \not = x'_{i+1}$, then $\lbrace x'_i, x'_{i+1} \rbrace
= \lbrace x_i, x_{i+1} \rbrace$. Therefore, for each integer $j$
with $x'_1 \leq j < x'_1 + k \vert q \vert e(K,D)$, the pair
$\lbrace j, j+1 \rbrace$ occurs as $\lbrace x'_i, x'_{i+1} \rbrace$
for some arc $\tilde C^1_i$. $\square$

The difficulty behind Theorem 5.1 is that there exist several
different normal curves $N$ in $\partial {\cal H}^0$ which have
zero combinatorial area. These are listed below (where we view each
curve $N$ as lying in $S^2$ and intersecting $G(D)$, since ${\cal F}$ is
dual to $G(D)$).

\vskip 12pt
\centerline{\psfig{figure=zeroarea.ps}}
\vskip 12pt
\centerline{Figure 23.}

We denote the normal curves $N$ in Figure 23 as $(k_1, k_2)$-curves,
where $k_1 = \vert N \cap {\cal F}^1 \vert$,
and $k_2 = \vert N \cap \partial M \vert$.
If $k_1 = 0$ and $k_2 = 2$, there are two possible
types of curve $N$, one of which is a boundary bigon.
We refer to the other type of curve as a $(0,2)$-curve.
We say that $C^0_i$ is a {\sl skirting arc} if it runs
between components of ${\cal F}^0 \cap \partial M$ which
are adjacent round $\partial {\cal H}^0 \cap \partial M$
(see Figure 24).

\vskip 12pt
\centerline{\psfig{figure=skirtarc.ps}}
\vskip 12pt
\centerline{Figure 24.}

\noindent {\bf Lemma 5.6.} {\sl If $C^0_i$ is a skirting arc,
then $H^0_i$ either has positive combinatorial area or
is a boundary bigon.}

\noindent {\sl Proof.} Suppose that, on the contrary,
$H^0_i$ has zero combinatorial area, but is not a boundary bigon.
By Lemma 4.2, $\partial H^0_i$ is normal, and hence is 
a $(2,1)$-curve or a $(0,2)$-curve.
This implies that each time $\partial H^0_i$ hits a vertex of $G(D)$, 
it crosses directly over that
vertex. This contradicts our assumption that $C^0_i$
is a skirting arc. $\square$

\noindent {\bf Lemma 5.7.} {\sl Suppose that $C^0_i$ and $C^0_{i+1}$ 
are both non-skirting arcs, and that $C^1_i$ runs along a 1-handle of
$\partial M$ corresponding to a non-bigon edge of $G(D)$. Then
at least one of $\partial H^0_i$ and $\partial H^0_{i+1}$ has positive combinatorial area.}

\noindent {\sl Proof.} Suppose that, on the contrary, both
$\partial H^0_i$ and $\partial H^0_{i+1}$ have zero combinatorial area. 
Then each is either a $(0,2)$-curve or a $(2,1)$-curve. View $\partial H^0_i$ and
$\partial H^0_{i+1}$ as curves in $S^2$ which intersect $G(D)$. Let
$A_i$ (respectively, $A_{i+1}$) be the sub-arc of
$\partial H^0_i$ (respectively, $\partial H^0_{i+1}$) corresponding to
the component of $\partial H^0_i \cap \partial H^1_i$ touching
$C^0_i \cap C^1_i$ (respectively, the component of $\partial H^0_{i+1}
\cap \partial H^1_i$ touching $C^0_{i+1} \cap C^1_i$).
Then $A_i$ and $A_{i+1}$ both lie
in the same region $R$ of $G(D)$, and they are matched onto each
other by a twist which is either clockwise or anti-clockwise,
depending on whether $R$ is coloured black or white.
There are a number of possibilities: 
\item{(i)} $A_i$ runs between successive vertices
of $R$;
\item{(ii)} $A_i$ runs between vertices of $R$ which are not successive;
\item{(iii)} $A_i$ runs between a vertex and an edge of $R$.

\noindent These cases are given in Figure 25.
In each case, it is clear that, since $D$ is prime and the regions of $G(D)$
admit a checkerboard colouring, the edge of
$G(D)$ corresponding to $C^1_i$ is a bigon edge. $\square$

\vskip 18pt
\centerline{\psfig{figure=primebig.ps}}
\vskip 12pt
\centerline{Figure 25.}

\noindent {\bf Lemma 5.8.} {\sl Suppose that precisely
one of $C_i^0$ and $C^0_{i+1}$ is a skirting arc, and that 
$C^1_i$ runs along a 1-handle of
$\partial M$ corresponding to a non-bigon edge of $G(D)$.
Suppose also that $x'_i \not=x'_{i+1}$. Then 
at least one of $H^0_i$ and $H^0_{i+1}$ has positive
combinatorial area.}

\noindent {\sl Proof.} Let $A_i$ and $A_{i+1}$ be as in the
proof of Lemma 5.7, and let $R$ be the region of $G(D)$ containing
them. By Lemma 5.6, if a skirting arc does not have
positive combinatorial area, then it is part of a boundary bigon.
Suppose therefore that one of $H^0_i$ and $H^0_{i+1}$ is
a boundary bigon. This implies that $A_i$ and $A_{i+1}$ both
run between successive vertices of $R$. Suppose (for definiteness)
that $A_i$ is parallel to the edge $E$ of $G(D)$ corresponding to 
$C^1_i$. If $\partial H^0_i$ is a $(0,2)$-curve, then as in case (i)
of the proof of Lemma 5.7, $E$ is
a bigon edge, contrary to assumption. Hence, $\partial H^0_i$
must be the boundary bigon, and $\partial H^0_{i+1}$
must be the $(0,2)$-curve. But, then, in fact
$x'_i = x'_{i+1}$, contrary to assumption. $\square$

A simple application of the formula for combinatorial area
gives the following result. The calculations are rather similar
to those in Proposition 4.10.

\noindent {\bf Lemma 5.9.} {\sl If $H^0_i$ has positive combinatorial
area, then $a(H^0_i) / \vert H^0_i \cap \partial M \vert \geq \pi/4$.}

The only type of normal curve $N$ for which 
$a(N) / \vert N \cap \partial M \vert = \pi/4$ is 
a $(1,2)$-curve given in Figure 26.

\vskip 12pt
\centerline{\psfig{figure=piover4.ps}}
\vskip 12pt
\centerline{Figure 26.}

\noindent {\sl Proof of Theorem 5.1.} 
Let $C$ be a curve representing a non-zero multiple $k$ of the
slope $p/q$ on $\partial {\cal N}(K)$, such that $C$
is disjoint from the 2-handles and respects the product
structure on the 1-handles. Define $C^0_1, \dots, C^0_n, C^1_1,
\dots, C^1_n$ as above and let ${\cal E} = H^0_1 \cup \dots \cup H^0_n
\cup H^1_1 \cup \dots \cup H^1_n$ be an inward extension of $C$.
Let $H$ and $H'$ be the 1-handles of $\partial {\cal N}(K)$ 
associated with a non-bigon edge of $K$. There are
$2 t(K,D)$ such non-bigon edges. By Lemma 5.5, there are
at least $k \vert q \vert$ arcs $C^1_i$ on $H \cup H'$
with $x'_i \not= x'_{i+1}$. Let $C^1_i$ be one such arc.
Then $H^0_i$ and $H^0_{i+1}$ cannot both be boundary bigons.
Lemmas 5.6, 5.7 and 5.8 imply
that at least one of $H^0_i$ and $H^0_{i+1}$ ($H^0_i$, say) has positive combinatorial
area. Lemma 5.9 gives that
$a(H^0_i) / \vert H^0_i \cap \partial M \vert \geq \pi/4$.
Hence, 
$$w({\cal E}) \geq (1/2) 2t(K,D) k \vert q \vert \pi/4 \geq t(K,D) \vert q \vert \pi/4.$$
As in the discussion after Proposition 4.10, the factor
(1/2) is necessary to avoid double counting. Since $C$ was
an arbitrary curve representing a non-zero multiple of the
slope $p/q$ on $\partial {\cal N}(K)$
and ${\cal E}$ was an arbitrary inward extension, we have that
$l(p/q) \geq \vert q \vert t(K,D) \pi/4$. $\square$

\vfill\eject
\centerline{\caps References}
\vskip 6pt

\item{1.} C. ADAMS, {\sl Augmented alternating link complements are
hyperbolic}, Low-dimensional Topology and Kleinian groups, London Math.
Soc. Lecture Note Series. {\bf 112}, Cambridge Univ. Press (1986)
\item{2.} I. AGOL, {\sl Ph.D. Thesis}, U.C. San Diego (1998)
\item{3.} I. AITCHISON, E. LUMSDEN and H. RUBINSTEIN, {\sl Cusp
structures of alternating links}, Invent. Math. {\bf 109} (1992)
473-494
\item{4.} S. BLEILER and C. HODGSON, {\sl Spherical space forms
and Dehn filling}, Topology {\bf 35} (1996) 809-833.
\item{5.} C. CAO and R. MEYERHOFF, {\sl The orientable cusped hyperbolic
manifolds of minimal volume}, Preprint.
\item{6.} C. DELMAN, {\sl Essential laminations and surgery on 2-bridge
knots}, Topology Appl. {\bf 63} (1995) 201-221.
\item{7.} D. B. A. EPSTEIN and R. C. PENNER, {\sl Euclidean decomposition
of non-compact hyperbolic manifolds}, J. Diff. Geom. {\bf 27} (1988) 
67-80.
\item{8.} D. GABAI, {\sl The ubiquitous nature of quasi-minimal
semi-Euclidean laminations in 3-manifolds}, Preprint.
\item{9.} D. GABAI and W. KAZEZ, {\sl Group negative curvature
for 3-manifolds with genuine laminations}, Geom. Topol. {\bf 2}
(1998) 65-77.
\item{10.} D. GABAI and U. OERTEL, {\sl Essential laminations in
3-manifolds}, Ann. Math. {\bf 130} (1989) 41-73.
\item{11.} M. GROMOV, {\sl Hyperbolic Groups}, Essays in Group Theory,
MSRI Publ. {\bf 8}, Springer (1987) 75-263.
\item{12.} J. HEMPEL, {\sl 3-Manifolds}, Ann. of Math. Studies, No. 86,
Princeton Univ. Press, Princeton, N. J. (1976)
\item{13.} C. HODGSON and S. KERCKHOFF, To appear.
\item{14.} W. JACO, {\sl Lectures on Three-Manifold Topology},  Regional
Conference Series in Mathematics, No. 43, Providence 1980, A. M. S.
\item{15.} W. JACO and U. OERTEL, {\sl An algorithm to decide if a 3-manifold
is Haken}, Topology {\bf 23} (1984) 195-209.
\item{16.} S. KOJIMA, {\sl Polyhedral decomposition of hyperbolic
3-manifolds with totally geodesic boundary}, Aspects of Low-dimensional
Manifolds, Adv. Stud. Pure Math. {\bf 20}, Kinokuniya, Tokyo (1992)
93-112.
\item{17.} W. MENASCO, {\sl Polyhedra representation of link complements},
Low- \hfill\break
dimensional Topology, Contemp. Math {\bf 20}, Amer. Math. Soc.
(1983)
\item{18.} W. MENASCO, {\sl Closed incompressible surfaces in alternating
knot and link complements}, Topology {\bf 23} (1984) 37-44.
\item{19.} W. MENASCO and M. THISTLETHWAITE, {\sl Surfaces with boundary
in alternating knot exteriors}, J. Reine Angew. Math. {\bf 426} (1992)
47-65.
\item{20.} L. MOSHER, {\sl Geometry of cubulated manifolds}, Topology {\bf
34} (1995) 789-814.
\item{21.} I. RIVIN, {\sl On geometry of convex ideal polyhedra in
hyperbolic 3-space}, Topology {\bf 32} (1993) 87-92.
\item{22.} R. ROBERTS, {\sl Constructing taut foliations}, Comment. Math.
Helv. {\bf 70} (1995) 516-545.
\item{23.} W. THURSTON, {\sl The geometry and topology of 
three-manifolds}, Princeton (1980)
\item{24.} W. THURSTON, {\sl Three-dimensional manifolds, Kleinian groups and
hyperbolic geometry}, Bull. Amer. Math. Soc. {\bf 6} (1982) 357-381.

\vskip 12pt
\+ DPMMS \cr
\+ University of Cambridge \cr
\+ 16 Mill Lane \cr
\+ Cambridge CB2 1SB \cr
\+ England. \cr

\end